\documentclass[12pt]{article}
\usepackage{}
\usepackage{comment}
\usepackage{epsfig}
\usepackage{latexsym}
\usepackage{caption}
\usepackage{subcaption}
\usepackage{amsfonts}

\usepackage{amssymb}
\usepackage{mathrsfs}
\usepackage{amsmath}
\usepackage{enumerate}
\usepackage{graphicx}
\usepackage{MnSymbol}
\usepackage{float}
\usepackage{pict2e}
\usepackage{enumerate}
\usepackage{graphics}
\usepackage{MnSymbol}
\usepackage{tikz}
\usepackage{cite}
\usepackage{ulem}
\numberwithin{equation}{section}

\usepackage{soul,color} 
\definecolor{lightgreen}{RGB}{153,255,153}

\definecolor{brightred}{RGB}{255,80,102}
 
\definecolor{lightpurple}{RGB}{229,204,255}

\soulregister\cite7
\soulregister\ref7
\soulregister\pageref7

\usepackage[colorinlistoftodos,prependcaption,textsize=scriptsize,color=olive!70, textwidth=30mm,]{todonotes}

\usepackage[pdfauthor={derajan},pdftitle={How to do this},pdfstartview=XYZ,bookmarks=true,
colorlinks=true,linkcolor=blue,urlcolor=blue,citecolor=blue,pdftex,bookmarks=true,linktocpage=true,  hyperindex=true]{hyperref}

\usepackage{color}

\makeatletter

\newcommand{\Rmnum}[1]{\expandafter\@slowromancap\romannumeral #1@}

\makeatother

\begin{document}

\newtheorem{theorem}{Theorem}[section]
\newtheorem{observation}[theorem]{Observation}
\newtheorem{corollary}[theorem]{Corollary}
\newtheorem{algorithm}[theorem]{Algorithm}
\newtheorem{definition}[theorem]{Definition}
\newtheorem{guess}[theorem]{Conjecture}
\newtheorem{claim}{Claim}[section]
\newtheorem{problem}[theorem]{Problem}
\newtheorem{question}[theorem]{Question}
\newtheorem{lemma}[theorem]{Lemma}
\newtheorem{proposition}[theorem]{Proposition}
\newtheorem{fact}[theorem]{Fact}


\makeatletter
  \newcommand\figcaption{\def\@captype{figure}\caption}
  \newcommand\tabcaption{\def\@captype{table}\caption}
\makeatother

\newtheorem{acknowledgement}[theorem]{Acknowledgement}

\newtheorem{axiom}[theorem]{Axiom}
\newtheorem{case}[theorem]{Case}
\newtheorem{conclusion}[theorem]{Conclusion}
\newtheorem{condition}[theorem]{Condition}
\newtheorem{conjecture}[theorem]{Conjecture}
\newtheorem{criterion}[theorem]{Criterion}
\newtheorem{example}[theorem]{Example}
\newtheorem{exercise}[theorem]{Exercise}
\newtheorem{notation}{Notation}
\newtheorem{solution}[theorem]{Solution}
\newtheorem{summary}[theorem]{Summary}

\newenvironment{proof}{\noindent {\bf
Proof.}}{\rule{3mm}{3mm}\par\medskip}
\newcommand{\remark}{\medskip\par\noindent {\bf Remark.~~}}
\newcommand{\pp}{{\it p.}}
\newcommand{\de}{\em}
\newcommand{\mad}{\rm mad}
\newcommand{\qf}{Q({\cal F},s)}
\newcommand{\qff}{Q({\cal F}',s)}
\newcommand{\qfff}{Q({\cal F}'',s)}
\newcommand{\f}{{\cal F}}
\newcommand{\ff}{{\cal F}'}
\newcommand{\fff}{{\cal F}''}
\newcommand{\fs}{{\cal F},s}

\newcommand{\wrt}{with respect to }
\newcommand{\q}{\uppercase\expandafter{\romannumeral1}}
\newcommand{\qq}{\uppercase\expandafter{\romannumeral2}}
\newcommand{\qqq}{\uppercase\expandafter{\romannumeral3}}
\newcommand{\qqqq}{\uppercase\expandafter{\romannumeral4}}
\newcommand{\qqqqq}{\uppercase\expandafter{\romannumeral5}}
\newcommand{\qqqqqq}{\uppercase\expandafter{\romannumeral6}}

\newcommand{\qed}{\hfill\rule{0.5em}{0.809em}}

\newcommand{\var}{\vartriangle}

\title{{\large \bf A survey on Hedetniemi's conjecture}}

\author{ Xuding Zhu\thanks{School  of Mathematical Sciences, Zhejiang Normal University,  China.  E-mail: xdzhu@zjnu.edu.cn. Grant Number: NSFC 12371359.  }}

\maketitle

\begin{abstract}
 In 1966, Hedetniemi conjectured that for any positive integer $n$ and graphs $G$ and $H$, if neither   $G$ nor $H$ is $n$-colourable, then $G \times H$ is   not $n$-colourable.   This conjecture has received significant attention over the past half century, and was 
disproved by Shitov in 2019. Shitov's proof shows that Hedetniemi's conjecture fails for sufficiently large $n$. Shortly after Shitov's result, smaller counterexamples were found in a series of papers, and it is now  known that Hedetniemi's conjecture fails for all $n \ge 4$, and  holds for $n \le 3$. Hedetniemi's conjecture has inspired extensive research, and many related problems remain open. This paper surveys the results and problems associated with the conjecture, and explains the ideas used in finding counterexamples.  
\end{abstract}

\section{Introduction}
\label{sec:1}

Products  of graphs  are basic and fundamental constructions in graph theory,  giving rise to important graph classes and deep structural problems. It is natural that properties and parameters of the product graph are closely related to or determined by those of the factor graphs. 
Various graph products and graph invariants  have been extensively studied   in the literature, among which Hedetniemi's conjecture specifically concerns  the categorical product   and the  chromatic number of graphs. 

The \emph{ categorical product}    $G \times H$ of  graphs $G$ and $H$  has vertex set $V(G) \times V(H)$, in which    $(x, y)  (x',y') \in E(G \times H)$   if and only if $xx' \in E(G)$ and $yy' \in E(H)$.  (We write  $uv \in E(G)$ or $u \sim_G v$   to mean that $u$ and $v$ are adjacent in $G$. When the graph $G$ is clear from the context, we may simply write $u \sim v$ for $u \sim_G v$).  This product is also known by other names  in the literature, such as  the \emph{tensor product} and the 
\emph{direct product}.

A {\em proper $n$-colouring} of a graph $G$ is a mapping $f: V(G) \to \{1,2,\ldots, n\}$ such that $f(x) \ne f(y)$ for any edge $xy$ of $G$. The {\em chromatic number} of $G$ is the minimum $n$ such that $G$ admits a proper $n$-colouring.
Given a proper $n$-colouring $\phi$ of $G$, one can define   a proper $n$-colouring $\Phi$  of $G \times H$ by setting $\Phi(x,y) = \phi(x)$. This implies that $\chi(G \times H) \le   \chi(G)$.  Similarly,  
$\chi(G \times H) \le \chi(H)$, and thus $\chi(G \times H) \le \min \{\chi(G), \chi(H)\}$.
In 1966, 
Hedetniemi conjectured  that  equality always holds \cite{Hed1966}.

\begin{conjecture} 
	\label{Conj-Hede} For all graphs $G$ and $H$,  $\chi(G \times H) = \min \{\chi(G), \chi(H)\}$.
\end{conjecture}

  Reference \cite{Hed1966} is a Ph.D.thesis that is not easily accessible.  The  conjecture  first appeared in  a journal article of Burr, P. Erd\H{o}s  and  Lov\'{a}sz \cite{BEL1976} in 1976. Ne\v{s}et\v{r}il and Pultr     \cite{NP1978}   attributed the conjecture to Lov\'{a}sz.  Duffus, Sands and Woodrow \cite{DSW1985} called it   ``Hedetniemi's conjecture". With many partial results being proved and the general case seemingly intractable, the conjecture has received significant   attention since the 1980's  \cite{Kla1996, Sau2001,Tar2008,Zhu1998}, and was finally refuted by Shitov \cite{Shi2019} in 2019.

In this paper, we survey  the  concepts, methods, results, and problems related to this conjecture.

\section{Homomorphism-monotone graph invariants}

\begin{definition}
	Assume $G$ and $H$ are graphs. 
    A {\em homomorphism} from $G$ to $H$ is a mapping $f: V(G) \to V(H)$ such that for every edge   $xy$ of $G$, $f(x)f(y)$ is an edge of $H$. We say $G$ is {\em homomorphic } to $H$, written as $G \to H$, if there is a homomorphism from $G$ to $H$.
\end{definition} 

 Two graphs $G$ and $H$ are {\em homomorphically equivalent}, written as $G \leftrightarrow H$, if $G \to H$ and $H \to G$. We denote by $\mathcal{G}$ the set of finite graphs. 
 Then $\leftrightarrow$ is an equivalence relation
 on $\mathcal{G}$. We denote by $\mathcal{G}/\leftrightarrow$
 the equivalence classes of finite graphs with respect to the equivalence relation $\leftrightarrow$. For $G \in \mathcal{G}$, let $[G]$ denote the equivalence class containing $G$.
 Then the homomorphism relation $\to$ is a partial order on $\mathcal{G}/\leftrightarrow$. 

A homomorphism from a graph $G$ to $K_n$ is equivalent to a proper $n$-colouring of $G$.
 Thus,  
graph
homomorphism is a generalization of graph colouring, and a homomorphism from $G$ to
$H$ is also called an $H$-colouring of $G$ \cite{HNbook}.

Note that for any graphs $G$ and $H$,  $G \times H \to G$ and $G \times H \to H$, with the projections be the homomorphisms.  

\begin{definition}
\label{def-hm}
    A graph invariant $\rho$ is {\em homomorphism-monotone}   if $G \to H$ implies that $\rho(G) \le \rho(H)$.
\end{definition}

 The chromatic number of graphs is an example of homomorphism-monotone invariants.

 Assume $\rho$ is a homomorphism-monotone graph invariant. As $G \times H \to G$, we have $\rho(G \times H) \le \rho(G)$. Similarly, $\rho(G \times H) \le \rho(H)$, and hence

 $$\rho(G \times H) \le \min \{\rho(G), \rho(H)\}.$$

For each homomorphism-monotone graph invariant $\rho$, it is natural to ask if the following $\rho$-version  Hedetniemi's type equality holds:

\begin{question}
	\label{Conj-rho-Hede} Is it true that for all graphs $G$ and $H$,  $\rho(G \times H) = \min \{\rho(G), \rho(H)\}$?
\end{question}

Motivated by Hedetniemi's conjecture, for many homomorphism-monotone graph invariants, the questions  have been studied.  
For  a few homomorphism-monotone graph invariants, the Hedetniemi-type equalities were proved and some of the equalities have important consequences. These include the fractional chromatic number, the strict vector chromatic number, the vector chromatic number and the coindex of the box complex of graphs. 
For some other homomorphism-monotone graph invariants, the questions have been studied and the problems remain open. These include the  Shannon OR-capacity, and the index of the box complex of graphs, quantum chromatic number and local chromatic number. The most important
result is that the equality does not hold in general for the chromatic number, as well as for the circular chromatic number.
 
  \subsection{Fractional chromatic number}
  
  There are a few equivalent definitions of the fractional chromatic number $\chi_f(G)$ of a graph $G$. One definition is through a linear program. We denote by $\mathcal{I}(G)$ the set of independent sets of $G$. A {\em fractional colouring }  of $G$ is a mapping $\phi: \mathcal{I}(G) \to [0,1]$ such that   for each vertex $v$ of $G$, $\sum_{v \in I, I \in \mathcal{I}(G)}\phi(I) \ge 1$. The {\em weight} of $\phi$ is $w(\phi) = \sum_{I \in \mathcal{I}(G)} \phi(I)$. The {\em fractional chromatic number} $\chi_f(G)$ of $G$ is defined as 
$$\chi_f(G) = \inf \{w(\phi): \phi \text{ is a fractional colouring of } G\}.$$

The fractional chromatic number of categorical product of graphs was first studied by this authro in 2002 \cite{Zhu2002}, where the focus is  whether the fractional version Hedetniemi's conjecture is true, and some special cases were verified.

 In 2005,   Tardif \cite{Tar2005a} proved that for any graphs $G$ and $H$, $\chi_f(G \times H) \ge \frac 14 \min\{\chi_f(G), \chi_f(H)\}$. In 2012, Zhang \cite{Zhang2012} studied the independent set of the product of vertex-transitive graphs, and the result  implies that for vertex transitive graphs $G$ and $H$, $\chi_f(G \times H) = \min \{\chi_f(G), \chi_f(H)\}$. The fractional version of Hedetniemi's conjecture was fully proved  by this author in \cite{Zhu2011}.     

 \begin{theorem} 
 	\label{thm-zhufrac} For any graphs $G$ and $H$,
 	$$\chi_f(G \times H) =  \min \{\chi_f(G), \chi_f(H)\}.$$
 \end{theorem}
 
 A
 feature of fractional chromatic number
 that plays a crucial role in the proof in \cite{Zhu2011} is the duality of linear programming. 
To prove that $\chi_f(G \times H) \ge \min\{\chi_f(G), \chi_f(H)\}$, we study the dual of the above linear program, which defines the {\em fractional clique number} of $G$. A {\em fractional clique} of $G$ is a map $f: V(G) \to [0,1]$
such that for any independent set $I$ of $G$, $\sum_{x \in I}f(x)
\le 1$. The \emph{total weight} of $f$ is $w(f)=f(V(G)) = \sum_{v \in V(G)}f(v)$. The {\em fractional clique number} $\omega_f(G)$ of $G$ is
the maximum total weight  of
a fractional clique $f$ of $G$.  
By the duality of linear programming, we
know that for any graph $G$, $$\chi_f(G) = \omega_f(G).$$  So Theorem \ref{thm-zhufrac} is equivalent to  
the following theorem.

 \begin{theorem} 
 	\label{thm-zhufrac2} For any graphs $G$ and $H$,
 	$$\omega_f(G \times H) =  \min \{\omega_f(G), \omega_f(H)\}.$$
 \end{theorem}

 The advantage of using the dual form is that one can guess a weight distribution to be a fractional clique of $G \times H$ that is of the correct total weight.  Given a maximum fractional clique of $G$ and a maximum fractional clique of $H$, there is a natural weight distribution as a candidate for a maximum fractional clique of $G \times H$. Such a weight distribution  was already given in \cite{Zhu2002}. 
 
Assume $g$ is a  maximum fractional clique    of $G$ and $h$ is a maximum
 fractional clique  of $H$.  Let $\phi: V(G \times H) \to [0,1]$ be defined as follows:
 $$\phi(x,y)=\frac{g(x)h(y)}{\max\{\omega_f(G), \omega_f(H)\}}.$$ 
 Obviously,
 $$w(\phi) =   \frac{\omega_f(G)\omega_f(H)}{\max\{\omega_f(G), \omega_f(H)\}} =
 \min\{\omega_f(G), \omega_f(H)\}.$$  So 
 $w(\phi) \ge \omega_f(G \times H) \ge \min \{\omega_f(G), \omega_f(H)\}.$
 
 It was conjectured in \cite{Zhu2002} that $\phi$ is a fractional clique of $G \times H$. To prove this conjecture, it   amounts to showing that   for any independent set $U$ of $G \times H$, 
 $$\sum_{(x,y) \in U} g(x)h(y) \le \max\{\omega_f(G), \omega_f(H)\}.$$
 This was proved 10 years later in \cite{Zhu2011}.

\subsection{Vector chromatic number, strict vector chromatic number and Shannon OR-capacity}

 The  vector chromatic number and strict vector chromatic number of a graph $G$ were first introduced by Karger,   Motwani  and   Sudan   in \cite{KMS1998}. 

 \begin{definition}
     \label{def-vecchi}
     Assume $t \ge 1$ is a real number. A {\em  vector $t$-colouring} (respectively, a {\rm strict vector $t$-colouring}) of $G$ is a mapping $\phi$ that assigns to each vertex $v$ a real vector $\phi(v)$ of length $t-1$ such that  for each edge $uv$ of $G$, $\phi(u)^T \phi(v) \le -1$   (respectively,  $\phi(u)^T \phi(v) = -1$).    The {\em   vector chromatic number } $\chi_v(G)$ (respectively, the {\em strict vector chromatic number} $\chi_{sv}(G)$) of $G$ is the infimum  $t$ such that $G$ admits a  vector $t$-colouring (respectively, a strict vector $t$-colouring). 
 \end{definition}

It was observed  in \cite{KMS1998} that the strict vector chromatic number of a graph $G$ is equal to the Lov\'{a}sz $\theta$-function of the complement $\bar{G}$  of $G$, i.e., 
 $\chi_{sv}(G)  = \theta(\bar{G})$. As noted in \cite{GRRSV2020}, the vector chromatic number of $G$ is equal to  the $\theta'$-function of the complement of $G$ introduced in \cite{Schrijver1979}.

It follows from the definition that $\chi_v(G) \le \chi_{sv}(G)$.   For the complete graph $K_n$, mapping the vertices of $K_n$ to the vertices of a regular $(n-1)$-dimensional simplex is a strict vector $n$-colouring. This is an optimal vector colouring of $K_n$. Hence $\chi_v(K_n)=\chi_{sv}(K_n)=n$. As $G$ is $n$-colourable if and only if $G \to K_n$, this implies that $\chi_{sv}(G) \le \chi(G)$ for any graph $G$.  It is also known that for any graph $G$, $\chi_{sv}(G) \le \chi_f(G)$.
 Therefore  for any graph $G$, $$\omega(G) \le \chi_v(G) \le \chi_{sv}(G) \le \chi_f(G) \le \chi(G).$$

 It is easy to see that both $\chi_v$ and $\chi_{sv}$ are homomorphism-monotone. 
The vector chromatic number and strict vector chromatic number of the categorical product 
of graphs were first studied by Godsil,   Roberson, \v{S}\'{a}mal, and  Severini in \cite{GRSS2016}, and they proved  that the Hedetniemi-type equality holds for the strict vector chromatic number.

\begin{theorem}
    \label{thm-chisv}
    For any graphs $G$ and $H$, 
$$\chi_{sv}(G \times H) = \min \{\chi_{sv}(G), \chi_{sv}(H)\}.$$ 
\end{theorem}

It was conjectured in \cite{GRSS2016} that the Hedetniemi-type equality holds for the vector chromatic number as well. This conjecture was confirmed by Godsil,   Roberson,   Rooney,   \v{S}\'{a}mal  and   Varvitsiotis  in \cite{GRRSV2020}.

 \begin{theorem}
    \label{thm-chiv}
    For any graphs $G$ and $H$, 
$$\chi_{v}(G \times H) = \min \{\chi_{v}(G), \chi_{v}(H)\}.$$ 
\end{theorem}

The vector chromatic number and the strict vector chromatic number of a graph can be formulated as semidefinite programs, and the proofs of Theorems \ref{thm-chisv} and \ref{thm-chiv}  rely on the duality property of  the corresponding semidefinite programs.

The Hedetniemi-type equality for another homomorphism-monotone graph invariant, the quantum chromatic number $\chi_q(G)$ of a graph $G$, was also studied in \cite{GRSS2016}. The  quantum chromatic number of a graph $G$ is defined through a graph homomorphism game. We refer the reader to \cite{GRSS2016} for the definition.  
It was proved in \cite{GRSS2016} that  the Hedetniemi-type equality $\chi_q(G \times H) = \min\{\chi_q(G), \chi_q(H)\}$ holds for  some special graphs $G$ and $H$, and it remains an open   question whether it holds for all graphs $G$ and $H$.

 Let the OR-product of two graphs $G$ and $H$, denoted by $G \cdot H$, be the graph with vertex set $V(G) \times V(H)$, in which $(x,y)(x',y')$ is an edge if either $xx' \in E(G)$ or $yy' \in E(H)$. Denote by $G^t$ the $t$-fold OR-product of $G$ and by   $\omega(G)$   the clique number of $G$.   
 The Shannon OR-capacity $C_{OR}(G)$ is defined as 
 $$C_{OR}(G) = \lim_{t \to \infty}  \sqrt[t]{\omega(G^t)}. $$

 The Shannon OR-capacity $C_{OR}$ is a homomophism-monotone graph invariant. Simonyi \cite{Sim2021} asked whether the Hedetniemi-type equality $C_{OR}(G \times H) = \min \{C_{OR}(G), C_{OR}(H)\}$ holds for all graphs $G$ and $H$. 
 
 There is not much support for a positive answer to this question. On the other hand, it was proved in \cite{Sim2021}  that if this equality does not hold for $C_{OR}$, then it also does not hold for   some other `nicer'  graph invariant.  We denote by   $G \oplus H$ the {\em join} of   $G$ and $H$, which  is the graph  obtained from the disjoint union of $G$ and $H$ by adding all edges between $V(G)$ and $V(H)$. 
 Given  a collection $\mathcal{S}$ of graphs that is closed under the join and the OR-product and that contains $K_1$, the {\em asymptotic spectrum} $Y(\mathcal{S})$ of $\mathcal{S}$ is the set of all maps $\phi: S \to \mathbb{R}_{\ge 0}$ for which the following hold:
 \begin{enumerate}
     \item $\phi(K_1)=1$.
     \item $\phi(G \oplus H) = \phi(G)+\phi(H)$.
     \item $\phi(G \cdot H) = \phi(G) \phi(H)$.
     \item If $G \to H$, then $\phi(G) \le \phi (H)$.
 \end{enumerate}
 A list of known elements in $Y(\mathcal{S})$ was given in \cite{Zui2019}. These include the fractional chromatic number $\chi_f(G)$ and the strict vector chromatic number $\chi_{sv}(G)$.
 On the other hand, $\chi(G)$ is not an element of $Y(\mathcal{S})$ as (3) is not satisfied by $\chi$.
 One may wonder if a graph invariant is more likely to satisfy the Hedetniemi's equality if it is an element of $Y(\mathcal{S})$.
 
 The Shannon OR-capacity  $C_{OR} $ is not an element in $ Y(\mathcal{S})$. It  was proved in  \cite{Hae1979} that (3) does not hold for $C_{OR}$, and proved in \cite{Alon1998}  that  (2) does not hold for $C_{OR}$. However, 
 it was proved in \cite{Zui2019} that for any $G \in \mathcal{S}$, $C_{OR}(G) = \min_{\phi \in Y(\mathcal{S})} \phi(G)$. 
 Using this result  in \cite{Zui2019}, Simonyi \cite{Sim2021} proved that if $G$ and $H$ are graphs for which 
  $C_{OR}(G \times ) < \min\{C_{OR}(G), C_{OR}(H)\}$, then there exists $\phi \in Y(\mathcal{S})$ such that $\phi(G \times H) < \min\{\phi(G), \phi(H)\}$.

\subsection{Index and coindex} 
 
 The topological method, which originated from  Lov\'{a}sz's celebrated proof of the Kneser conjecture   \cite{Lov1978},  is a widely used tool in the study of graph colouring problems. The method involves associating  a topological space to a graph, where certain topological invariants of the associated space provide lower bounds for the chromatic number of the graph. The {\em box complex} ${\rm Box}(G)$ of a graph $G$ is a simplicial complex with vertex set $V({\rm Box}(G)) = V(G) \times \{+,-\}$. A  set $(A \times \{+\}) \cup (B \times \{-\})$ forms a simplex if $A,B \subseteq V(G)$ and $ab \in E(G)$ for all $a \in A, b \in B$. The space ${\rm Box}(G)$ is a {\em $Z_2$-space},  equipped with   a $Z_2$-action $\nu: V({\rm Box}(G)) \to V({\rm Box}(G))$ defined as $\nu((v,+)) = (v,-)$ and $\nu((v,-))=(v,+)$. 

 For two $Z_2$-spaces $X,Y$, with $Z_2$-actions $\nu_X, \nu_Y$, respectively, a {\em $Z_2$-map} from $X$ to $Y$ is a continuous map $\phi: X \to Y$ satisfying $\phi(\nu_X(x)) = \nu_Y(\phi(x))$ for all $x \in X$. We write $X \to_{Z_2} Y$ if such a map exists, and write $X \not\to_{Z_2} Y$ otherwise.

 Denote by $S^d$  the $d$-dimensional sphere, with $Z_2$-action $\nu$ defined as $\nu(x) = -x$ for each point $x \in S^d$. It is known and easy to verify that for any positive integer $n$,   ${\rm Box}(K_n) \simeq_{Z_2} S^{n-2}$.    The famous 
 Borsuk-Ulam theorem  states that 
	for any positive integer $d$, $S^d \not\to_{Z_2} S^{d-1}$. This motivates the following definitions.

 \begin{definition}
     \label{def-index}
     For a $Z_2$-space $X$, the  {\em index } of $X$, denoted by  ${\rm ind}(X)$, is the smallest integer $n$ such that $X \to_{Z_2} S^n$. The  {\em coindex } of $X$, denoted by  ${\rm coind}(X)$, is the largest integer $n$ such that $S^n \to_{Z_2} X$.
 \end{definition}

 The index ${\rm ind}(Box(G))$ and the coindex ${\rm coind}(Box(G))$ of the box complex $Box(G)$ of a graph $G$ provide lower bounds for the chromatic number $\chi(G)$. We refer to \cite{MR1988723} for detailed introduction of ${\rm ind}(Box(G))$ and the coindex ${\rm coind}(Box(G))$, and other  topological spaces associated with graphs and related topological invariants. 
    For our purpose, both parameters ${\rm ind}(Box(G))$ and   ${\rm coind}(Box(G))$ can also be defined via graph homomorphisms.   

  Given a positive integer $k$ and $0 < \epsilon < 2$,  the {\em  Borsuk graph} $B_{k,\epsilon}$ is the graph   whose vertices are the points of the $k$-dimensional unit sphere $S^{k}$,  in which $xy$ is an edge  if the distance between $x$ and $y$ (in $\mathbb{R}^{k+1}$) is at least $2-\epsilon$. 
  
   It was shown in \cite{ST2006} that  for any graph $G$, $${\rm coind}(G) = \max \{k: \exists \epsilon > 0, B_{k,\epsilon} \to G\}.$$

 Alternately, ${\rm coind}(G)$ can also be defined by using homomorphisms from generalized Mycielski graphs. For a positive integer $n$, let $P_n$ be the path with vertices $0,1,\ldots, n$. For a graph $G$, the {\em $n$th generalized mycielskian} of $G$, denoted by $M_n(G)$, if the graph obtained from $G \times P_n$ by adding edges $\{(x,0)(y,0): xy \in E(G)\}$ and identifying all vertices $\{(x,n): x \in V(G)\}$ into a single vertex. 
For $k \ge 2$, the family $\mathcal{K}_k$ of generalized Mycielski graphs is defined recursively as follows: $\mathcal{K}_2 =\{K_2\}$. For $k \ge 3$, $\mathcal{K}_k = \{M_n(G): G \in \mathcal{K}_{k-1}, n \in \{1,2,\ldots \}\}$. It was shown in \cite{STW2017} that for any graph $G$, 
 $${\rm coind}(G) = \max \{k: \exists G' \in \mathcal{K}_k, G' \to G\}.$$

The  two formulas can be used as   alternate definitions of ${\rm coind}(Box(G))$. 
It follows from these definitions that 
${\rm coind}(Box(G)) $ is a homomorphism monotone invariant.

 Borsuk-Ulam Theorem states that if sphere $S^k$ is covered by $k+1$ open sets, then one of the open sets contains a pair of anti-podal points, i.e., two points in $S^k$ of distance 2. It follows from this theorem that $\chi(B_{k, \epsilon}) \ge k+2$. As remarked by Lov\'{a}sz \cite{MR708798},  Borsuk-Ulam Theorem is indeed equivalent to the statement that $\chi(B_{ k, \epsilon}) \ge k+2$ for any $0 < \epsilon < 2$.   This lower bound is sharp if $\epsilon$ is small enough \cite{ST2006,MR1988723}.   Thus for any graph $G$, we have 
 $$\chi(G) \ge {\rm coind}(Box(G)) +2.$$

To define ${\rm ind}(Box(G))$ via graph homomorphisms, we need to introduce an important graph operation.

For a graph $K$ and a positive integer $d$, let $\Omega_{2d+1}(K)$ be the graph whose vertices are $(d+1)$-tuples  $(A_0, A_1, \ldots, A_d)$ such that 
\begin{enumerate}
	\item $A_i \subseteq V(K)$ for $i=0,1,\ldots, d$.
	\item $|A_0|=1$   and $A_1 \ne \emptyset$.
	\item $A_i \subseteq A_{i+2}$ for $i=0,1,\ldots, d-2$.
	\item $A_d$ and $A_{d-1}$ are fully adjacent.
\end{enumerate}
Two tuples $(A_0, A_1, \ldots, A_d)$ and $(B_0, B_1, \ldots, B_d)$ are adjacent if $A_i \subset B_{i+1}$ and $B_i \subseteq A_{i+1}$ for $i=0,1,\ldots, d-1$, and $A_d$ and $B_d$ are fully adjacent.

The operator $\Omega_{2d+1}$ is an important graph functor, which plays an important role in the construction of counterexamples to Hedetniemi's conjecture.
It is difficult to draw or visualize the graph $\Omega_{2d+1}(G)$ even for small graphs $G$. One simple example is that $\Omega_{2d+1}(C_{2k+1}) = C_{(2d+1)(2k+1)}$. So $\Omega_{2d+1}$ acts
on cycles as  a subdivision of the edges. For an arbitrary graph $G$,  $\Omega_{2d+1}$ behaves much like a subdivision in the topological sense on the box complex. They preserve the homotopy type and refine the geometric structure. Continuous maps between box complexes can be approximated with graph homomorphisms from the refinement $\Omega_{2d+1}(G)$ of $G$. 

The graph  $\Omega_{2d+1}(K_n)$ is also denoted by $W(d+1,n)$ in the literature, and can be constructed as follows (see \cite{Wro2020}): 

The vertices of $\Omega_{2d+1}(K_n)$   are tuples $(x_1,x_2, \ldots, x_n) \in \{0,1,\ldots, d+1\}^n$ such that 
There exists exactly one $i \in [n]$ for which $x_i=0$, there exists some $i \in [n]$ for which $x_i=1$,  and $(x_1,x_2,\ldots, x_n)$ is adjacent to $(y_1, y_2, \ldots, y_n)$ if and only if for all $i \in [n]$, either $|x_i-y_i|=1$ or $x_i=y_i=d+1$.

The following two  theorems  were proved by Wrochna \cite{Wro2019}.

\begin{theorem}
    \label{thm-index}
    There exists a $Z_2$-map 
  from ${\rm Box}(G)$ 
 to ${\rm Box}(H)$
 if and only if for some positive integer $k$, there is a homomorphism from $\Omega_{2k+1}(G)$ to $H$. 
\end{theorem}

\begin{theorem}
    \label{thm-wro-omega-equiv}
    $|{\rm Box}(G)|$ and $|{\rm Box}(\Omega_{2d+1}(G))|$ are $Z_2$-homotopy equivalent for any graph $G$ with no loops and for any positive integer $d$.  
\end{theorem}   

The following property of $Box(G)$ can be used as an alternate definition of ${\rm ind}(Box(G))$:
$${\rm ind}({\rm Box}(G)) = \min\{n: \text{ there exists a positive integer $k$ such that } \Omega_{2k+1}(G) \to K_{n+2}\}.$$

It follows from the definitions that ${\rm ind}(Box(G))$ is a homomorphism monotone graph invariant, and for any graph $G$, 
$$\omega(G) \le {\rm coind}({\rm Box}(G))+2 \le {\rm ind}({\rm Box}(G)) +2\le \chi(G).$$

The lower bounds ${\rm coind}({\rm Box}(G))+2 $ and/or  $ {\rm ind}({\rm Box}(G)) +2$ for $\chi(G)$ are referred to as topological lower bounds. For some families of graphs whose fractional chromatic number and chromatic number are far apart, these  bounds can be significantly better than the other (non-topological) bounds. 
In particular, these bounds are sharp for Kneser graphs.    
It is easy to see that $B_{ k, \epsilon} \to B_{ k, \epsilon'}$ if  $\epsilon \le \epsilon'$. 
Using this observation, the following result was proved by Simonyi and Zsb\'{a}n \cite{SZ2010}.

\begin{theorem}
    \label{thm-SZ}
    For any graphs $G$ and $H$, 
    $$  {\rm coind}({\rm Box}(G \times H)) = \min\{{\rm coind}({\rm Box}(G )),{\rm coind}({\rm Box}(  H)).$$
\end{theorem}

The following result is a consequence of Theorem \ref{thm-SZ}, and was proved earlier by Hell  \cite{Hell1979}. 

\begin{theorem}
    \label{thm-Hell}
    If   $\chi(G) = {\rm coind}(Box(G)) +2$ and $\chi(H) = {\rm coind}(Box(H)) +2$, then $\chi(G \times H)=\min\{\chi(G), \chi(H)\}$.
\end{theorem}

 Wrochna \cite{Wro2019} studied the question whether the Hedetniemi-type equality holds for the index of graphs.
It was proved in \cite{Ziv2005} and \cite{Cso2007} that any $\mathbb{Z}_2$-space is $\mathbb{Z}_2$-equivalent to the ${\rm Box}(G)$ for some graph $G$, and it
 was proved in \cite{Cso2008} that $|{\rm Box}(G)| \times |{\rm Box}(H)| \simeq_{\mathbb{Z}_2} |{\rm Box}(G \times H)|$ for any graphs $G$ and $H$, where for 
  two $\mathbb{Z}_2$-spaces $(X, \nu_X)$ and $(Y, \nu_Y)$, $X \times Y$ is a $\mathbb{Z}_2$-space with $\mathbb{Z}_2$ action $\nu$ defined as $\nu(x,y) = (\nu_X(x), \nu_Y(y))$. It follows from these results that the Hedetniemi-type equality for the indices of box complex of graphs is equivalent to the same equality for the indices of all $\mathbb{Z}_2$-spaces. The following conjecture was proposed by Wrochna \cite{Wro2019}:

\begin{conjecture}
    \label{conj-wrochna}
    For all $Z_2$-spaces $X$ and $Y$, ${\rm ind}(X \times Y) = \min\{ {\rm ind}(X), {\rm ind}(Y)\}$.
\end{conjecture}

By the results mentioned above, we may restrict the $Z_2$ spaces $X$ and $Y$ in Conjecture \ref{conj-wrochna}  to box complexes of graphs, and hence this conjecture  is equivalent to the following conjecture \cite{Wro2019}:
     
     \begin{conjecture}
     \label{conj-top-equiv}
         For any graph $G$ and $H$, if $\chi(G \times H) \le n$, then for some integer $d$, $\Omega_{2d+1}(G)$ or $\Omega_{2d+1}(H)$ is $n$-colourable.
     \end{conjecture}

     We may compare Conjecture \ref{conj-top-equiv} with Hedetniemi's conjecture.  Under the same condition that $\chi(G \times H) \le n$, Hedetniemi's conjecture asserts that $G$ or $H$ is $n$-colourable, which is significantly stronger.


More generally, 
a $\mathbb{Z}_2$-space $Z$ is called   {\em multiplicative} if $X \times Y \to_{\mathbb{Z}_2} Z$ implies that $X \to_{\mathbb{Z}_2} Z$ or $Y \to_{\mathbb{Z}_2} Z$. Conjecture \ref{conj-wrochna} says that $S^d$ is multiplicative for any non-negative integer $d$. We say a graph $K$ is {\em multiplicative } if for any graphs $G$ and $H$, $G \times H \to K$ implies that $G \to K$ or $H \to K$.   Hedetniemi conjecture says that $K_n$ is multiplicative for all positive integer $n$.   
The following result was proved independently by Matsushita \cite{Matsushita} and Wrochna \cite{Wro2019}.

\begin{theorem}
    \label{thm-multi-Z2space}
    If $K$ is a multiplicative graph, then $|{\rm Box}(K)|$ is a multiplicative $\mathbb{Z}_2$-space.
\end{theorem}

It follows from Theorem \ref{thm-multi-Z2space} that Hedetniemi's conjecture is stronger than Conjecture \ref{conj-wrochna}. 
As $K_2$ and $K_3$ are multiplicative, it follows from Theorem \ref{thm-multi-Z2space} that $S^0$ and $S^1$ are multiplicative. 
For $n \ge 4$, we know that $K_n$ is not multiplicative. On the other hand, it remains open whether $S^d$ is multiplicative for all $d \ge 2$. 

Conjecture \ref{conj-wrochna} was studied by Daneshpajouh,   Karasev, and   Volovikov in \cite{DRV2023}, where it was proved that the conjecture holds if one of the factor spaces 
is a  sphere.  Moreover, the conjectured  equality holds with index replaced by   the
homological index. 
      
\subsection{Circular chromatic number} 

Given positive integers $p \ge 2q$, a {\em  $(p,q)$-colouring} of a graph $G$ is a mapping $f: V(G) \to \{0,1,\ldots, p-1\}$ such that for every edge $xy$ of $G$, 
$$q \le |f(x)-f(y)| \le p-q.$$
The {\em circular chromatic number} of $G$ is defined as 
$$\chi_c(G) = \inf\left\{ \frac pq: G \text{ admits a $(p,q)$-colouring} \right\}.$$
The circular chromatic number of a graph $G$ is a refinement of its chromatic number introduced by Vince in \cite{Vin1988}. The term ``circular chromatic number" was coined in \cite{ZhuSurvey2001}.  Observe that a $(p,1)$-colouring is the same as a proper $p$-colouring. Hence $\chi_c(G) \le \chi(G)$ for every graph $G$. On the other hand, it was shown in \cite{Vin1988} that $\chi_c(G) > \chi(G)-1$. Hence $\chi(G) = \lceil \chi_c(G) \rceil$. 
It is also known and easy to see that $\chi_c(G) \ge \chi_f(G)$ for every graph $G$. 

As a strengthening of Hedetniemi's conjecture, this author conjectured in \cite{Zhu1992} that for any graphs $G$ and $H$, $\chi_c(G \times H) = \min\{\chi_c(G), \chi_c(H)\}$. As Hedetniemi's conjecture is false, so is this one. Nevertheless, the result that Hedetniemi's conjecture holds for $n \le 3$, i.e.,  $\chi(G \times H)\le 3 \rightarrow \chi(G) \le 3 $ or $\chi(H) \le 3$  is extended to $r < 4$, i.e., $\chi_c(G \times H) \le r < 4$ implies that $\chi_c(G) \le r$ or $\chi_c(H) \le r$. This is related to the concept of multiplicative graphs, which we shall discuss more in Section 4.

\subsection{Chromatic Ramsey number}

Hedetniemi's conjecture was independently proposed in \cite{BEL1976}, where the conjecture was inspired by a problem on  the  chromatic  Ramsey number of graphs.
 
 Suppose $G,H,F$ are graphs. We write $F \longrightarrow   (G,H)$ if for any
 colouring of the edges of $F$ with colours red and blue, there is
 either a red copy of $G$ (i.e., $G$ is a subgraph (not necessarily induced subgraph) of the red graph)
 or a blue copy of $H$. The \emph{Ramsey number} of $G,H$ can be defined as 
 $$R(G, H) = \min\{|V(F)|: F
 \longrightarrow  (G,H)\}.$$
 By replacing the parameter $|V(F)|$ with other graph parameters, one defines a variety of Ramsey parameters. For example, 
 $$R_E(G, H) = \min\{|E(F)|: F
 \longrightarrow  (G,H)\}$$ 
 is the \emph{size Ramsey number} of $G,H$. 
 The \emph{chromatic Ramsey number} 
 $R_{\chi}(G,H)$ is defined as $$R_{\chi}(G, H) = \min\{\chi(F): F
 \longrightarrow   (G,H)\}.$$ We write $R_{\chi}(G)$ for $R_{\chi}(G,G)$. It was proved in \cite{BEL1976}  that 
 $$R_{\chi}(G) = \min \{n: \forall c: E(K_n) \to [2], \exists i \in [2], G \to   K_n[c^{-1}(i)]\}.$$
 Here $K_n[c^{-1}(i)]$ denotes the subgraph of $K_n$ induced by edges coloured by colour $i$.
 In other words, $R_{\chi}(G)$ is the minimum $n$ such that every 2-edge colouring of $K_n$ contains a monochromatic subgraph which is a homomorphic image of $G$.
 
 If $\chi(F) \le (n-1)^2$, and $c$ is an $(n-1)^2$-colouring of $F$
 with colours $\{(i,j): 1\le i,j\le n-1\}$, then we can colour the
 edges of $F$ as follows: If $e=xy$, $c(x)=(i,j)$ and $c(y)=(i',j')$
 then we colour $e$ red if $i=i'$ and colour $e$ blue otherwise. It
 is easy to see that both the red graph and the blue graph are
 $(n-1)$-chromatic. This shows that if  $\chi(G)=n$ then $R_{\chi}(G)
 \ge (n-1)^2+1$.  Burr, Erd\H{o}s and Lov\'{a}sz \cite{BEL1976} 
 conjectured that this lower bound is sharp. Let $M(n) = \min\{R_{\chi}(G): \chi(G) = n\}$.

 \begin{conjecture}
     \label{conj-bel}
      For any positive integer $n$, $M(n) = (n-1)^2+1$.
 \end{conjecture}
 
 A tentative proof of this conjecture was sketched in \cite{BEL1976}:

 Let $m=(n-1)^2+1$. Let $c_1, c_2, \ldots, c_p$ be the set of all 2-edge colourings of $K_m$. For $i \in [p]$, let $G_{i,1}, G_{i,2}$ be the two monochromatic subgraphs of $K_m$ under colouring $c_i$. The following claim is easily verified:

 \begin{claim}
     \label{clm-km}
     For 2-edge colouring of $K_m$, one of the monochromatic subgraphs $G$ has $\chi(G) > n-1$.  
 \end{claim}

 Therefore for each $i \in [p]$, there  exists $j_i \in [2]$ such that $\chi(G_{i,j_i}) > n-1$. 
 Let 
 	$$G=G_{1, j_1} \times G_{2, j_2} \times \ldots \times G_{p, j_p}.$$
 	
  If $\chi(G) \ge n$, then for  any 2-edge colouring $c_i$ of $K_m$, $G_{i, j_i}$ is a monochromatic subgraph, which is a homomorphic image of $G$ (the $i$th projection map is a homomorphism from $G$ to $G_{i, j_i}$). Therefore $R_{\chi}(G) \le m = (n-1)^2+1$, and hence  $R_{\chi}(G)  = (n-1)^2+1$. 

  The problem  is  why should $\chi(G_{1, j_1} \times G_{2, j_2} \times \ldots \times G_{p, j_p}) \ge n$. As each $G_{i,j_i}$ has chromatic number at least $n$, this inequality would follow from the seemingly  natural general identity $\chi(G \times H) = \min\{\chi(G), \chi(H)\}$, which they proposed as a conjecture   \cite{BEL1976}. 

  Now this conjecture is refuted. However, this conjecture  is stronger than the required inequality $\chi(G_{1, j_1} \times G_{2, j_2} \times \ldots \times G_{p, j_p}) \ge n$, as each $G_{i,j_i}$ are ``special" graphs, which satisfies not only  the inequality $\chi(G_{i,j_i}) \ge n$, but also some other inequalities. 
  
  Tardif \cite{Tar2024} noted that for the proof to work, it suffices to have a graph invariant $\rho(G)$ which is a lower bound for $\chi(G)$ for which (i)  Claim \ref{clm-km} holds, i.e., for any 2-edge colouring of $K_m$, one of the monochromatic subgraphs $G$ has $\rho(G) > n-1$,   
   and  (ii) the Hedetniemi-type equality holds.
 
   If $K_m$ is 2-edge coloured by red and blue and the red subgraph $G_R$ has fractional chromatic number at most $n-1$, then $G_R$ has clique number at most $n-1$. This implies that the blue subgraph $G_B$ has independence number $\alpha(G_B) \le n-1$, and hence has fractional chromatic number at least $\frac{m}{\alpha(G_B)} > n-1$. Using this observation and Theorem \ref{thm-zhufrac}, the following theorem was proved in \cite{Zhu2011}.

    \begin{theorem} 
      \label{thm-zhuchromaticramsey}
      For any positive integer $n$, $M(n)=(n-1)^2+1$. 
  \end{theorem} 
  
  It was proved by Lov\'{a}sz in \cite{Lov1979} that for any graph $G$, $\theta(G) \theta(\bar{G}) \ge |V(G)|$, where $\theta(\bar{G}) = \chi_{sv}(G)$ is the Lov\'{a}sz theta number of $\bar{G}$ and  the strict vector chromatic number of $G$. Therefore if $K_m$ is 2-edge coloured, then one of the monochromatic subgraphs $G$ has $\chi_{sv}(G) > n-1$.  
  Using this fact and Theorem \ref{thm-chisv}, an alternate proof of Theorem \ref{thm-zhuchromaticramsey} was given by Tardif in \cite{Tar2024}.

\section{Counterexamples to Hedetniemi's conjecture}

One of the key tools used in the study of Hedetniemi's conjecture  is the concept of  exponential graphs.  

\begin{definition}
	Assume $K$ and $G$ are graphs. The exponential graph $K^G$ has vertex set 
	$$V(K^G) = \{f: V(G) \to V(K) \}$$
	and edge set $$E(K^G)= \{fg: \forall xy \in E(G), f(x)g(y) \in E(K)\}.$$
\end{definition}

In particular, $ K^G$ has a loop if and only if $G \to K$.  

The following is a fundamental property of the exponential graph $K^G$:

\begin{proposition}
	\label{prop-1}
	For any graphs    $G,H,K$, $G \times H \to K$ if and only if $H \to K^G$.
\end{proposition}
The mapping $\phi: V(K^G \times G) \to V(K)$ defined as $\phi(f,x) = f(x)$ is a homomorphism from $K^G \times G \to K$. 
On the other hand, if $\phi$ is a homomorphism from $G \times H$ to $ K$, then $\psi: V(H) \to V(K^G)$ defined as $\psi(u) (x) = \phi(x,u)$ is a homomorphism from $H$ to $K^G$.

In other words, in the homomorphism order, $K^G$ is a maximum graph among all graphs $H$ for which $G \times H \to K$.

For the study of Hedetniemi's conjecture, we are particularly interested in the exponential graph $K_n^G$. We use $[n]=\{1,2,\ldots, n\}$ as the vertex set of $K_n$. The vertices of $K_n^G$ are all the $n$-colourings (not necessarily proper) of $G$, in which two $n$-colourings $f$ and $g$ are adjacent if for every edge $xy$ of $G$, $f(x) \ne g(y)$. 

As a particular case of Proposition \ref{prop-1}, 
we know that for any graphs $G$ and $H$, $G \times H$ is $n$-colourable if and only if 
	$H \to K_n^G$.  
 
 For a positive integer $n$, let $H(n)$ be the following statement:

 \bigskip
\noindent
$H(n)$:  
\emph{ If neither $G$ nor $H$ is $n$-colourable, then $G \times H$ is not $n$-colourable.  }
\bigskip
 
 Hedetniemi's conjecture says that $H(n)$ holds for all positive integer $n$. 
 If $\chi(G) > n$ and $\chi(K_n^G)> n$, then $H(n)$ fails because $ G \times K_n^G $ is $n$-colourable. If $\chi(G)> n$ and $\chi(K_n^G) \le n$, then for any graph $H$ with $\chi(H) > n$, $H$ is not homomorphic to $K_n^G$ and hence $\chi(G \times H) > n$. Thus  Hedetniemi's conjecture is equivalent to the following statement:

\bigskip
\noindent
$H'(n)$:  
\emph{ For any  graph $G$, if $\chi(G) > n$ then $\chi(K_n^G) \le n$.}
\bigskip

Shitov \cite{Shi2019} proved  that $H'(n)$ fails if $n$ is sufficiently large. No explicit bound for $n$ was given in \cite{Shi2019}, and a careful analysis of the proof shows that   $n \ge 2 \cdot 83^2 \cdot 3^{82}$ suffices. Shortly after Shitov's work, smaller  counterexamples were found in a series of papers, and $H'(n)$ was shown to be false for $n \ge 125$ in  \cite{Zhu2020}, for $n \ge 13$ in  \cite{Tar2022c}, for $n \ge 5$ in  \cite{Wro2020}, and finally for $n \ge 4$ in \cite{Tar2023}. 

The results in these paper got not only 
progressively stronger, but also the proofs got progressively simpler (at least to read), culminating in the really surprisingly short proof by Tardif.

It was proved by El-Zahar and Sauer \cite{ES1985} that $H'(n)$ holds for $n \le 3$. Thus the problem of whether $H'(n)$ holds for a given positive integer $n$ is completely solved. 

To find a counterexample for $H'(n)$, it amounts to finding a graph $G$ for which  $\chi(G) > n$ and $\chi(K_n^G) >n$. 

The counterexample constructed by Shitov uses the lexicographic product of graphs. The {\em lexicographic product} of two graphs $F$ and $H$, denoted by $F[H]$, is obtained from $F$ by replacing each vertex of $F$ with a copy of $H$.    In particular,  $F[K_q]$ is obtained from $F$ by   blowing up each vertex of $F$ into a clique of size $q$. Thus the vertex set of
$F[K_q]$ is  $\{(x,i): x \in V(F), i \in [q]\}$,   in which $(x,i)$ and $ (y,j)$ are adjacent   if and only if 
either $x y \in E(F)$ or $x=y$ and $i \ne j$.

The   {\em girth} of $F$ is the length of a shortest  cycle in $F$, and the {\em odd girth} of $F$ is the length of a shortest odd cycle in $F$. Shitov  proved the following theorem in \cite{Shi2019}:

\begin{theorem}\label{thm-shitov}
    Assume $q$ is sufficiently large, $F$ is a graph of odd girth $7$ with $\chi_f(F) > 3+2/q$. Let $G=F[K_q]$ and $n=3q+2$. Then $\chi(G)> n$ and $\chi(K_n^G) > n$.
\end{theorem}

The existence of graphs with arbitrary large girth and arbitrary large fractional chromatic number was proved  in \cite{Erd1959}.  A computer search done by Exoo (see \cite{Zhu2021}) shows that the circulant graph $  {\rm Cay} (\mathbb{Z}_{83}, \{ \pm 2,  \pm 22, \pm 25, \pm 34 \})$
has odd girth $7$ and has fractional chromatic number  $\chi_f(F) = \frac{83}{27} > 3.07$ (and hence $\chi_f(F)> 3+2/q$ for $q \ge 29$).  

We call a vertex $\phi$ of $K_n^G$   a   mapping, as $\phi$ is just a mapping from $V(G)$ to $[n]$. Given a mapping $\phi$, let $Im(\phi)=\{\phi(x): x \in V(G)\}$ be the range of $\phi$. Note that if $\phi$ and $\phi'$ are mappings with $Im(\phi) \cap Im(\phi') = \emptyset$, then $\phi$ is adjacent to $\phi'$ in $K_n^G$. When $G=F[K_q]$, we say a mapping $\phi$ is {\em simple}  if $\phi(x,i)=\phi(x,j)$ for any $x \in V(F)$ and $i,j \in [q]$.
For a simple mapping, we may write $\phi(x)$ for $\phi(x,i)$.

It is well-known that for any graph $F$, 
$$\chi_f(F) = \inf\left\{ \frac{\chi(F[K_q])}{q}: q \in \mathbb{N}\right\}.$$
Therefore $\chi(G) \ge \chi_f(F)q > 3q+2=n$. To prove Theorem \ref{thm-shitov}, 
the difficult part is to prove that $\chi(K_n^G) > n$.

Shitov's proof and the proofs in the subsequent papers have a common feature: A small  subgraph $H$ of $K_n^G$ is shown to be  not $n$-colourable. 

This is rather natural. The graph $K_n^G$ is huge, but a large part of the graph is very sparse. In particular, there are many isolated vertices: Assume $\phi$ is a vertex of $K_n^G$. If there is a vertex $x$ of $G$ for which $\phi(N_G(x))=[n]$, then $\phi$ is an isolated vertex of $K_n^G$. On the other hand, $K_n^G$ contains an $n$-clique induced by the constant mappings $\{{\rm const}_i: i\in [n]\}$, where ${\rm const}(x)=i$ for all 
vertices $x$ of $G$. The subgraph $H$ of $K_n^G$ used in the proofs is induced by some vertices in the neighbourhood of this $n$-clique. Indeed, the important vertices of $H$ are those vertices $
\phi$ that are similar to constant mappings: $\phi(x)=c$ is constant $c$ for most vertices $x$ of $G$. Usually the image $Im(\phi)$ is small.

In the following, we describe some main ideas in these  proofs.   First consider Theorem  \ref{thm-shitov}.

Assume to the contrary that $K_n^G$ has a proper $n$-colouring $\Psi$. Since $\{ {\rm const}_i: i \in [n]\}$ induces a clique in $K_n^G$, we may assume $\Psi({\rm const}_i)=i$ for all $i\in [n]$.
Then for any $\phi \in V(K_n^G)$, $\Psi(\phi) \in Im(\phi)$, because if $ i \notin Im(\psi)$, then $\psi$ is adjacent to $    {\rm const}_i$, and hence $\Psi(\psi) \ne \Psi({\rm const}_i)=i$.

For each vertex $v$ of $F$, for any two distinct colours $b,t \in [n]$, let 
$\theta_{v,b,t} \in V(K_n^G)$ be  a simple mapping  defined as 
\[
\theta_{v,b,t}(x) = \begin{cases} b, & \text{ if $d_F(x,v)\ge 2$,} \cr
t, & \text{ if $d_F(x,v)\le 1$}.
\end{cases}
\]

As $Im(\theta_{v,b,t}) =\{b,t\} $, we know that $\Psi (\theta_{v,b,t}) \in \{b,t\}$. 
The colour $t$ is used in a small portion of $G$, i.e., vertices $(x,i)$, where $x$ is in the closed neighbourhood of $v$. Intuitively, $\theta_{v,b,t}$ is ``close" to the constant mapping ${\rm const}_b$.
So   $\Psi (\theta_{v,b,t}) $ is more likely to be $b$. 
This intuition is justified   in the following claim.

\begin{claim}
	\label{clm1}
	For each vertex $v$, there exists a colour $t=t(v) \in \{q+1, q+2, \ldots, 3q+2\}$ such that for each $b \in \{2q+1,2q+2,\ldots, 3q+2\} - \{ t\}$, 
	$\Psi(\theta_{v,b,t})=b$.
\end{claim}

For $v \in V(F)$ and $b \in \{2q+1,2q+2,\ldots, 3q+2\} - \{ t\}$,  let 
  $$I(v,b)=\{\phi \in V(K_n^G): \text{ $\phi$ is a simple mapping and}   \Psi(\phi)=b = \phi(v) \}.$$

The mapping $\theta_{v,b,t}$ is used to bound the size of $I(v,b)$.

Assume $t=t(v)$ and $b \in \{2q+1,2q+2,\ldots, 3q+2\} - \{ t\}$. As   $\Psi(\theta_{v,b,t})=b=\Psi(\phi)$,  any $\phi \in I(v,b)$ is not adjacent to $\theta_{v,b,t}$. This means that there exists $x \in V(F)$ and $y \in N_F[x]$ ($N_F[x]$ denotes the closed neighbourhood of $x$ in $F$), $\phi(x) = \theta_{v,b,t}(y)$. 

Observe that $x \ne v$,
  because $\phi(v)=b$, and for any $y \in N_G[v]$, $\theta_{v,b,t}(y) =t \ne b$. So there are $p-1$ choices for $x$. Once $x$ is chosen,  $\phi(x)$ has only two choices: $\phi(x)  \in  Im(\theta_{v,b,t}) = \{b,t\}$. To determine a mapping $\phi \in I(v,b)$, there is one choice for $\phi(v)$ (i.e., $\phi(v) = b$), there are $p-1$ choices for $x$ and two choices for $\phi(x)$, and for each vertex $u \in V(F)- \{v,x\}$, there are   $n$ choices for $\phi(u)$. This implies that $$|I(v,b)| \le 2(p-1)n^{p-2}.$$

Now we can bound $| \bigcup_{v \in V(F),  b \in \{2q+1,2q+2,\ldots, 3q+2\} - \{ t(v)\} }I(v,b)|$ as follows: 
there are $p$ choices for $v$ and $n/3$ choices for $b$, hence
$$\left| \bigcup_{v \in V(F),  b \in \{2q+1,2q+2,\ldots, 3q+2\} - \{ t(v)\} }I(v,b) \right| \le \sum_{v \in V(F),  b \in \{2q+1,2q+2,\ldots, 3q+2\} - \{ t(v)\} }|I(v,b)| \le 2(p-1)n^{p-2} \cdot p   n/3.$$
Let 
$$\mathcal{I} = \{\phi: \text{ $\phi$ is a simple mapping and $\phi(v) \in [n] - \{1,2,\ldots, 2q, t(v)\}$ for each vertex $v \in V(F)$} \}. $$
Then for each $\phi \in \mathcal{I}$,   $\Psi(\phi) = \phi(v) $ for some $v \in V(F)$ and $b = \phi(v) \in \{2q+1,2q+2,\ldots, 3q+2\} - \{ t(v)\}$. 
Therefore $\phi \in I(v,b)$ for some $v$ and $b$. Hence 
$$\mathcal{I} \subseteq \bigcup_{v \in V(F),  b \in \{2q+1,2q+2,\ldots, 3q+2\} - \{ t(v)\} }I(v,b). $$
As there
are at least $n/3$ choices for $\phi(v)$ for each vertex $v \in V(F)$, we conclude that $$|\mathcal{I} | \ge \left(\frac{n}{3}\right)^p.$$
This implies that $$\frac 23 p(p-1)n^{p-1} \ge \left(\frac{n}{3}\right)^p.$$
However, this is not true   when $n=3q+2 >  2p^23^{p-1}$ (i.e., $q \ge 2p^23^{p-2}$), and hence  Theorem \ref{thm-shitov} is proved.

A crucial point in the proof above is that   $I(v,b)$ is relatively small. This is due to the fact that   each  $\phi \in I(v,b)$ is coloured the same colour as $\theta_{v,b,t}$, and hence is not adjacent to $\theta_{v,b,t}$.    The final counting argument shows that when $n$ is sufficiently large,   it is impossible  all these $I(v,b)$ be small. In other words,  the following proposition holds:

\begin{proposition}
    \label{prop-2}
    If $q$ is large enough, then for some $v$ and $b$, 
there exists a simple mapping $\phi \in I(v,b)$
which is adjacent to $\theta_{v,b,t}$.  
\end{proposition}

In \cite{Shi2019}, the existence of the simple mapping $\phi$ was proved by a counting argument. In \cite{Zhu2021}, a simple mapping $
\phi$ is explicitly constructed, which is adjacent to all relevant vertices $\theta_{v,b,t}$. Indeed,  the mappings $\theta_{v,b,t}$ in \cite{Shi2019} are replaced by mappings $\theta_{v,t}$, defined as 
\[
\theta_{v,t} (x) = \begin{cases} v, & \text{ if $d_F(x,v)\ge 2$,} \cr
t, & \text{ if $d_F(x,v)\le 1$},
\end{cases}
\]
where the vertex $v$ also plays the role of the colour $b$. Let $\phi$ be the  simple mapping  defined by $\phi(v)=v$. Then $\phi$  is adjacent to   $\theta_{v,t}$  for all $v$ and all relevant $t$. This dramatically reduces the quantity  $q$ (and hence  $n$). The following result was proved in \cite{Zhu2021}.

\begin{theorem}
	\label{thm-Zhu}
	Let $F$ be a graph of odd girth at least $7$, with $p$ vertices.  Then for   $q \ge (p-1)/2$  and $n=3q+2$,  
	$\chi \left  ( K_n^{F[K_q]} \right ) > n$. 
\end{theorem}

 \begin{figure}[!htb]
 	\centering
 	\includegraphics[scale=0.65]{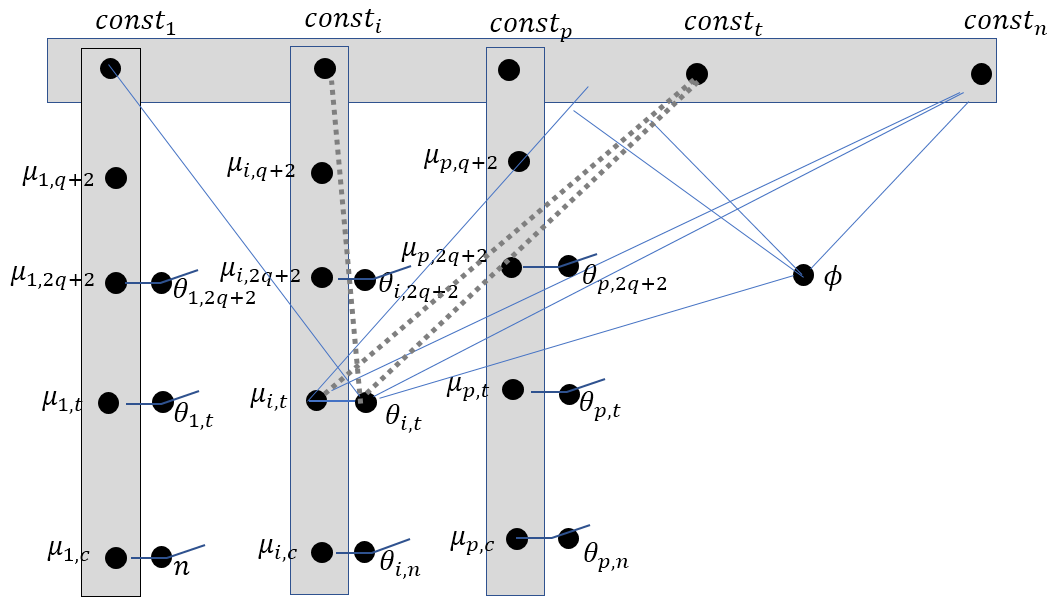}
 	\caption{A non-$n$-colourable subgraph $H$ of $K_n^G$ for Theorem  \ref{thm-Zhu} }\label{fig1}
 \end{figure}

In the proof of Theorem \ref{thm-Zhu}, a 
  non-$n$-colourable subgraph $H$ of $K_n^G$  is constructed explicitly, and is illustrated in Figure \ref{fig1}. Its  structure is quite simple. We omitted the definition of $\mu_{i,t}$. To have a full proof of Theorem \ref{thm-Zhu}, we do need all the definitions and need to check adjacency between these vertices. These checking are not difficult, but a little tedious. For the proof that $H$ is not $n$-colourable, all we need are the following:
  \begin{itemize}
      \item Each shaded rectangle is a clique. The top rectangle is the clique induced by the $n$ constant maps. 
      \item ${\rm Im}(\phi) = [p]$ and $\phi$ is adjacent to $\theta_{i,t}$ for all $i$ and $t$.
      \item $Im(\theta_{i,t}) = \{i,t\}$, and $\theta_{i,t}$ is adjacent to $\mu_{i,t}$.
      \item $Im(\mu_{i,t}) \subseteq [p] \cup \{t\}$.
  \end{itemize}
  Assume to the contrary that $\Psi$ is an $n$-colouring of $H$ with $\Phi(f) \in Im(f)$ for each mapping $f$. The mapping $\phi$ has image $Im(\phi) = [p]$, and hence $\Psi(\phi) = i$ for some $ i \in [p]$. The clique induced by $\{ \mu_{i,t}: t= q+2, q+3, \ldots, n\} \cup \{{\rm const}_i\}$ has size $n-q = 2q+2=p+1$. Hence   one of the mappings $\mu_{i,t}$ is coloured by $t$ for some $t \ge p+1$. Now $\theta_{i,t}$ is adjacent to both $\phi$ and $\mu_{i,t}$, but the two colours $i,t$ in $Im(\theta_{i,t})$ are used by $\phi$ and $\mu_{i,t}$, respectively. Therefore we obtain a contradiction.

Let $F$ be  the $83$-vertex graph  of odd girth $7$ with $\chi_f(F) > 3.07$ (hence  $p=83$), and let   $q=41$ and $n=3q+2 = 125$, 
we conclude  that $H'(n)$ fails for $n=125$.


It is possible that there exists a $p$-vertex graph $F$ of odd girth $7$   with $\chi_f(F) >  3+ 4/(p-1)$ for a smaller $p$. If this is the case, then   $H'(n)$ would  fail for a smaller $n$.  However, it is unlikely that $p$ can be significantly reduced. 

 Tardif realized that the number $p$ in Theorem \ref{thm-Zhu} need not be the number of vertices of $F$. Instead one can let $p$ be the number of colour classes of $F$ in a proper colouring of $F$, provided that for each colour class $X$,  $N^{=i}(X)$ is   an independent set for $i=1,2$. Here 
 $$N^{=i}(X)=\{y: \exists x \in X, \text{ and } F \text{ has an $x$-$y$-walk of length $i$} \}.$$
Indeed such a colouring is used in \cite{Zhu2021},   where all vertices of $F$ are coloured with distinct colours, and each colour class is a single vertex. The condition that    $N^{=i}(\{x\})$ is an independent set for $i=1,2$  is equivalent to say that $F$ has no triangle and no 5-cycle, i.e., $F$ has odd girth $7$.

  A $p$-colouring $\gamma$ of a graph $F$ is called a {\em $k$-wide $p$-colouring} of $F$ if   for each colour $i$, 
 $N^{=k'}(\gamma^{-1}(i))$ is an independent set of $F$ for all $k' < k$.

 Tardif proved the following result in \cite{Tar2022c}:
 
 \begin{theorem}
 	\label{thm-Tardif} 
 	Let $q, c, m$ be positive integers such that $n \ge \max\{m+q+1, 3q+2\}$. If $F$ is a graph which has a $3$-wide $m$-colouring, then $\chi(K_n^{F[K_q]}) > n$.
 \end{theorem} 

  The proof of Theorem \ref{thm-Tardif} is essentially  the same as the proof of  Theorem \ref{thm-Zhu}, except that in the definition of the mappings,  the distance to a single vertex is replaced by the distance to a colour class.

 To show that $H'(n)$ fails for $n$   using Theorem \ref{thm-Tardif}, one needs to find a graph $F$ that has a 3-wide $m$-colouring and  $\chi(F[K_q]) > n = \max\{m+q+1, 3q+2\}$.

 Indeed, in the proof of Theorem \ref{thm-Zhu}, such a graph $F$, a 3-wide colouring of $F$, and the  integers $q$ and $n$ satisfy this relation are used, where   $m = |V(F)|$ and the 3-wide $m$-colouring assigns each vertex  a distinct colour.  
 
 To apply Theorem \ref{thm-Tardif} to show that $H'(n)$ fails for a smaller $n$, one needs to find such a graph $F$  that has a 3-wide $m$-colouring for some small integer $m$, while ensuring that $F[K_q]$ has a sufficiently large chromatic number. The number of vertices of $F$ is irrelevant.

 For a graph $G$ and a positive integer $k$, let $\Gamma_{2k+1}(G)$ be the graph obtained from $G$ by adding edges $xy$ for each pair of vertices $x,y$ that are connected by a walk of length exactly $2k+1$. The graph $\Gamma_{2k+1}(G)$ is called the $(2k+1)$th power of $G$ and it is sometimes denoted by $G^{2k+1}$ (this is different from the graph obtained from $G$ by connecting every pair of vertices of distance at most $2k+1$, which is also usually  denoted by $G^{2k+1}$). A $(k+1)$-wide $p$-colouring of $G$ is the same as a $p$-colouring of $\Gamma_{2k+1}(G)$. 
 
 The problem whether there is an $n$-chromatic graph $G$ which has a $3$-wide $n$-colouring  was first studied by Gy\'{a}rf\'{a}s, Jensen and Stiebitz \cite{GJS2004}. 
 The problem of existence of $n$-chromatic graph $G$ admitting  a $d$-wide $n$-colouring  for arbitrary $d$ was studied by Baum and Stiebitz \cite{BS2005},  Hajiabolhassan 
\cite{Haj2009},  and Wrochna \cite{Wro2019}. 

Topological methods were used in answering these questions.
It follows from Theorem \ref{thm-wro-omega-equiv} that  for any positive integer $d$, $\chi(\Omega_{2d-1}(K_n))=n$, and it was   proved in \cite{Haj2009,Wro2019} that 
for any graph $G$, $G$ admits a $d$-wide $n$-colouring if and only if $G$ is homomorphic to $\Omega_{2d-1}(K_n)$.

Let $q=4, n=14, m=9$ and $F=\Omega_5(K_m)$. By Theorem \ref{thm-Tardif}, $\chi(K_n^{F[K_q]}) > n$.    Lov\'{a}sz \cite{Lov1978} proved that for any positive integers  $m \ge 2q$,  $\chi(K(m,q)) = m-2q+2.$ 
Thus $\chi(K(14,4))=8 < 9= \chi(F)$, and hence $F$ is not homomorphic to $K(14,4)$. Note that 
an $n$-colouring of $F[K_q]$ is equivalent to a homomorphism from $F$ to $K(n,q)$. 
Therefore $\chi(F[K_4]) > n$, and  $H'(n)$ fails for $n =14$. 

Tardif \cite{Tar2022c}   further proved that $H'(n)$ fails for $n = 13$, by adding another idea.
   
However, if one  attempts to apply Theorem \ref{thm-Tardif} to prove $H'(n)$ fails for $n=13$, one needs to choose $q=3$.
The problem with this choice is that the chromatic number of $(\Omega_{2d+1}(K_m))[K_q]$ is not large enough. 
 
 Gujgiczer and Simonyi \cite{GS2022} proved that $(\Omega_{2d+1}(K_m))[K_q]$ has chromatic number $m+2q-2$  for $d+1 \ge q \ge 2$. This implies that for $q=3, n=m+q+1$, and $F=\Omega_5(K_m)$ the required inequality $\chi(F[K_q]) > n$ is not satisfied.

Wrochna \cite{Wro2020} addressed this problem by swapping the operations $\Omega_{2d+1}$ and the lexicographic product. Instead of considering the graph $(\Omega_{2d+1}(K_m))[K_q]$, he considered the graph $\Omega_{2d+1}(K_m[K_q]) = \Omega_{2d+1}(K_{mq})$, and proved the following theorem in \cite{Wro2020}.

\begin{theorem} 
	\label{thm-Wrochna}
	Assume $m,q,n$ are positive integers such that $m \ge  q+1$ and $n \ge m+q+1$. If $G$ admits a 3-wide $mq$-colouring,  then $\chi(K_n^G) > n$. 
\end{theorem}
 
  Let $q=2, n=7,m=4$ and $G=\Omega_5(K_{n+1})$. Then the  conditions of Theorem \ref{thm-Wrochna} are satisfied. As $\chi(G)>n$, we conclude that $H'(n)$ fails for $n \ge 7$.   
 
 By considering $\Omega_{2d+1}(K_m)$ for large $d$, Wrochna \cite{Wro2020} improved the result above and showed that $H'(n)$ fails for $n \ge 5$.  
 Tardif further proved in \cite{Tar2022} that  for any $n \ge 5$, there exists an infinite family 
 $\{G_m: m \in \mathbb{N}\}$   of graphs such that $\chi(G_m) > n$ for all $m \in \mathbb{N}$  and $\chi(G_m \times G_{m'}) \le n$ for all $m \ne m'$. 
 
 Finally, by also  considering $\Omega_{2d+1}(K_m)$ for large $d$, Tardif \cite{Tar2023} showed that $H'(n)$ fails for $n \ge 4$. Indeed, Tardif proved the following   stronger result in \cite{Tar2023}.

 \begin{theorem} \label{thm-Tardif2}
     For any integers $n \ge 4$ and $m=2^k$,   $\chi(K_n^{\Omega_{4k+1}(K_m)}) > n$.
 \end{theorem}


To prove that $H(n)$ fails for $n$, one needs to find graphs $G$ and $H$ such that
$\chi(G), \chi(H) > n$ and $\chi(G \times H) \le n$. Theorem \ref{thm-Tardif2} implies that there are graphs $G$ and $H$ such that $\chi(G) > m, \chi(H) > n$ and $\chi(G \times H) \le n$, where $m$ can be  arbitrarily large. In this sense, Theorem \ref{thm-Tardif2} is stronger than just saying $H(n)$ fails for $n \ge 4$.

\begin{figure}[!htb]
 	\centering
 	\includegraphics[scale=0.65]{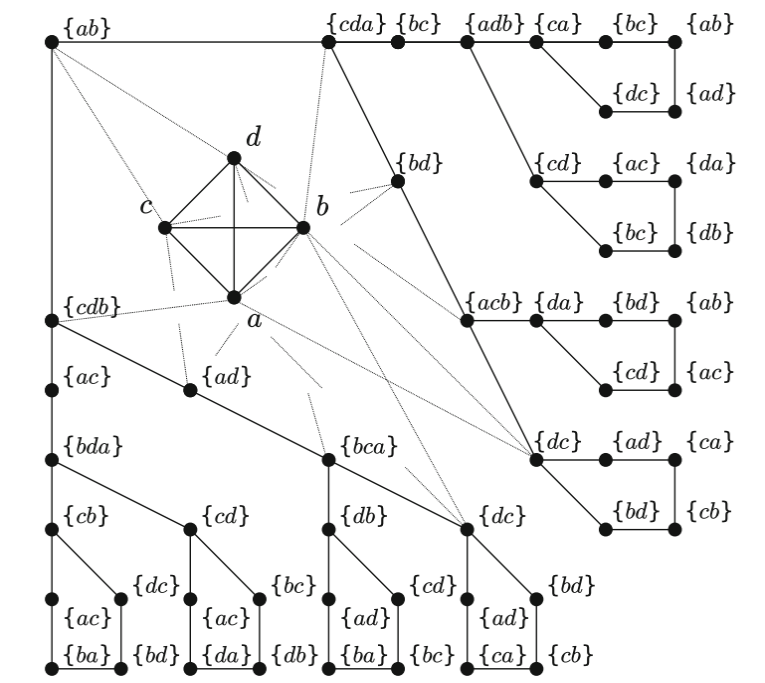}
 	\caption{A non-$4$-colourable subgraph $H$ of $K_4^{\Omega_{13}(K_8)}$, where a vertex labeled $\{a\}$ indicates the constant mapping ${\rm const}_a$, a vertex labeled $\{ab\}$ indicates  a mapping $\sigma_{a,b}^{\ell, S}$ for some $\ell \in \{0,1,\ldots, 6\}$ and $S \subseteq \{1,2,3,4\}$,  and a vertex labeled $\{abc\}$  indicates  a mapping $\tau_{a,b,}^{\ell, S, T}$ for some $\ell \in \{0,1,\ldots, 6\}$ and $S, T \subseteq \{1,2,3,4\}$. The labels of vertices also indicate allowed colors, as in list coloring of the graph.  }\label{fig2}
 \end{figure}

To prove Theorem \ref{thm-Tardif2},   a small subgraph $H$ of $K_n^{\Omega_{4k+1}(K_m)}$ (as illustrated in Figure \ref{fig2}) is shown to be not $n$-colourable.

We present the vertices in the subgraph $H$. It requires deep insight to find this subgraph. However, it is not difficult to check the adjacency among these vertices and to show that 
the subgraph is not $n$-colourable.

For a subset $S$ of  $[m]$,  colours $a,b \in [n]$ and $\ell \in \{0,1,\ldots, 2k\}$, let $\sigma_{a,b}^{\ell, S} \in V(K_n^{\Omega_{4k+1}(K_m)})$ defined as  
 \[
 \sigma_{a,b}^{\ell, S}(X_0,X_1,\ldots, X_{2k}) = \begin{cases} 
 a, &\text{ if $S \cap X_{\ell} \ne \emptyset$} \cr b, &\text{ otherwise.}
 \end{cases} 
 \]
 For disjoint subsets $R,T$ of $[m]$,   colours $a,b,c \in [n]$, and $\ell \in \{0,1,\ldots, 2k\}$, let $\tau_{a,b,c}^{\ell, R,T}  \in V(K_n^{\Omega_{4k+1}(K_m)})$ be   defined as 
 \[
 \tau_{a,b,c}^{\ell, R,T}(X_0,X_1,\ldots, X_{2k}) = \begin{cases} 
 a, &\text{ if $R \cap X_{\ell} \ne \emptyset$} \cr b, &\text{ if $R \cap X_{\ell} = \emptyset$ and $T \cap X_{\ell} \ne \emptyset$} \cr
 c, &\text{ otherwise.}
 \end{cases} 
 \]

 Recall that $\sigma_{a,b}^{\ell, S}$ and $\tau_{a,b,c}^{\ell, R,T}$   are vertices of $K_n^{\Omega_{4k+1}(K_m)}$, and  hence are functions from $V(\Omega_{4k+1}(K_m))$   to $[n]$.
 
 It is proved in \cite{Tar2023}  that for $\ell \in \{0,1,\ldots, 2k-1\}$,   distinct colours $a,b,c,d$,   disjoint subsets $R,T$ of $[m]$, and for   $S=R \cup T$,  the following adjacency hold: 
 $$\sigma_{a,b}^{\ell, S} \sim \sigma_{c,a}^{\ell, S} \ {\rm and } \ \sigma_{a,b}^{\ell, S} \sim \tau_{c,d,a}^{\ell+1,R,T} \ {\rm and } \  \tau_{a,b,c}^{\ell,R,T} \sim \sigma_{d,a}^{\ell+1, R}, \sigma_{d,b}^{\ell+1, T}.$$ 
 Additionally, for   $x \in [m]$, distinct colours $a,b,c \in [n]$, $\ell \in \{0,1,\ldots, 2k\}$, we have
 $$\sigma_{a,b}^{\ell, 
 \{x\}} \sim \sigma_{a,c}^{\ell, \{x\}}.$$ 
 
  To prove that $\chi(K_n^{\Omega_{4k+1}(K_m)}) > n$,   assume that 
   $\Psi$ is an $n$-colouring of $K_n^{\Omega_{4k+1}(K_m)}$,  with $\Psi(\phi) \in Im(\phi)$. It is proved (Lemma 8 of \cite{Tar2023}) by induction on $i=0,1,\ldots, k-1$, that 
 there exists a set $S$ of size $2^{k-i-1}$ and distinct colours $a,b$ such that $\Psi(\sigma_{a,b}^{2i,S})=a$. 
 Thus there exists $x \in V(K_m)$ and colours $a,b \in [n]$ such that   $\Psi(\sigma_{a,b}^{2k-2,\{x\}}) = a$. 
 Using the adjacency relations listed above, we conclude that $$ [\sigma_{b,c}^{2k,\{x\}}, \sigma_{c,a}^{2k-1,\{x\}}, \sigma_{a,b}^{2k-2,\{x\}}, \sigma_{d,a}^{2k-1,\{x\}}, \sigma_{b,d}^{2k,\{x\}}]$$
 forms a 5-cycle. 
Since $\Psi(\sigma_{a,b}^{2k-2,\{x\}})=a$ and 
$\Psi ( \sigma_{c,a}^{2k-1,\{x\}} ) \in \{c,a\}$, we conclude that $\Psi ( \sigma_{c,a}^{2k-1,\{x\}} ) = c$. Similarly, we deduce that $\Psi(\sigma_{b,c}^{2k,\{x\}})= \Psi(\sigma_{b,d}^{2k,\{x\}}) = b$. But $\sigma_{b,c}^{2k,\{x\}} \sim \sigma_{b,d}^{2k,\{x\}}$, a contradiction.

Note that  the vertices $\phi$ of $K_n^G$ that are used in the above argument have $|Im(\phi)| \le 3$.  In particular, all these vertices are in the closed neighbourhood of the $n$-clique induced by the constant maps. 
 

For the graphs $G$ and $H$ for which $\chi(G \times H) = n < \min\{\chi(G), \chi(H)\}$ constructed above, one of $G$ and $H$ has  clique number $n$. It is well-known \cite{NZ2004} that for any graph $G$ and any positive integer $g$, there is a graph $G'$ of girth at least $g$ such that $G' \to G$ and $\chi(G')=\chi(G)$. As observed by Tardif, using this result, one can obtain graphs $G'$ and $H'$ of arbitrarily large girth  such that $\chi(G' \times H') = n < \min\{\chi(G'), \chi(H')\}$.

Based on Theorem \ref{thm-zhufrac},  a weakening of Hedetniemi's conjecture was proposed in \cite{Zhu2011}:   if $\chi_f(G) > n $ and $\chi(H) > n$, then $\chi(G \times H) > n$. The counterexamples in Theorem \ref{thm-shitov}, Theorem \ref{thm-Zhu} and Theorem \ref{thm-Tardif} are all counterexamples to this weaker conjecture. 

Based   Theorem \ref{thm-Hell}, one may ask whether 
  $ {\rm coind}(G) +2 > n$ and $\chi(H) > n$ imply $\chi(G \times H) > n$. The counterexamples in Theorem \ref{thm-Wrochna} and Theorem  \ref{thm-Tardif2} give negative answer to this question.

\section{Multiplicative graphs and digraphs}

Recall that a graph $K$ is {  multiplicative} if for any two graphs $G$ and $H$, $G \not\to K$ and $H \not\to K$ implies that $G \times H \not\to K$.
  
It is easy to see that $K_1$ and $K_2$ are multiplicative. The celebrated result of
 El-Zahar and Sauer \cite{ES1985} states  that $K_3$ is multiplicative.

\begin{theorem}
    \label{thm-es}
    $K_3$ is multiplicative. 
\end{theorem}  

Now we know that $K_n$ is not multiplicative for all $n \ge 4$.  The problem of the multiplicativity of complete graphs is completely solved. 
However, the question of which other graphs  are multiplicative remains widely open.

The concept of  homomorphism and multiplicative graphs are naturally generalized to digraphs. For a digraph $G$, its vertex set and arc set are  denoted by $V(G)$ and $E(G)$, respectively. An arc (or a directed edge) of $G$ is denoted by an ordered pair $(x,y)$.  Digraphs in this paper have no parallel arcs, but may have pairs of opposite arcs. 

A homomorphism from a digraph $G$ to a digraph $H$ is a mapping $f:V(G) \to V(H)$ such that for any arc $(x,y)$ of $G$, $(f(x),f(y))$ is an arc of $H$.  The categorical product  $G \times H$ of digraphs $G$ and $H$ has vertex set $\{(x,y): x \in V(G), y \in V(H)\}$, in which $((x,y), (x',y'))$ is an arc if and only if $(x,x')$ is an arc of $G$ and $(y,y')$ is an arc of $H$.  A digraph $K$ is multiplicative if $G \not\to K$ and $H \not\to K$ implies that $G \times H \not\to K$.
We denote by $\mathcal{D}$ the family of digraphs. Similarly, two digraphs $G$ and $H$ are homomorphically equivalent, written as $G \leftrightarrow H$, if $G \to H$ and $H \to G$. We denote by $\mathcal{D}/\leftrightarrow$ the equivalence classes of digraphs with respect to the equivalence relation $\leftrightarrow$, and denote by $[G]$ the equivalence class in $\mathcal{D}/ \leftrightarrow$ that contains $G$. Homomoprhism relation is a partial order on $\mathcal{D}/\leftrightarrow$.

A digraph $K$ is   {\em symmetric} if $(x,y) \in E(K)$ implies that $(y,x) \in E(K)$. A symmetric digraph is equivalent to an undirected graph, where each pair of opposite arcs is equivalent to an undirected edge. Note that the product $G \times H$ of two symmetric digraphs is a symmetric digraph. In this sense,  the  family $\mathcal{G}$ of graphs is a subfamily of the family $\mathcal{D}$ of digraphs. The homomorphism and categorical product of undirected graphs are simply the restriction of homomorphism and categorical product of digraphs   to symmetric digraphs. 

Note that multiplicativity is sensitive to the family of digraphs under consideration.  
When we say $K_3$ is multiplicative in Theorem \ref{thm-es}, it means $K_3$ is multiplicative in the family of symmetric digraphs, or equivalently  in the family of undirected graphs. 
The precise statement is that ``$K_3$ is multiplicative in $\mathcal{G}$''. Indeed, there are digraphs $G$ and $H$ such that $G \times H \to K_3$ and $G \not\to K_3$ and $H \not\to K_3$. Thus $K_3$ is not multiplicative in $\mathcal{D}$. 

If $G \not\to H$ and $H \not\to G$, then $G \times H$ is a non-multiplicative graph. These are trivial non-multiplicative graphs. It seems more difficult to construct examples of multiplicative graphs, and also difficult to construct non-trivial non-multiplicative graphs. Theorem \ref{thm-es} 
was generalized by H\"aggkvist et al. \cite{HHMN1988} to odd cycles.

\begin{theorem}
    \label{thm-hhmn}
    Every odd cycle is multiplicative in $\mathcal{G}$.
\end{theorem}

 The proof of Theorem \ref{thm-hhmn} is an adaptation of the proof of Theorem \ref{thm-es}. For a long time, these were the only known multiplicative graphs. In \cite{Tar2005}, Tardif proved that the operation $\Omega_{2d+1}$ can be used to construct new multiplicative graphs from known multiplicative graphs.

It is more convenient to put the graph constructions in the setting of category theory. A {\em digraph functor } (respectively, a {\em  graph functor}) is a construction $\Gamma$ that constructs a  digraph $\Gamma(G)$ from a  digraph $G$ (respectively, a  graph $\Gamma(G)$ from a   graph $G$)  such that if there is a homomorphism from $G$ to $H$, then there is a homomorphism from $\Gamma(G)$ to $\Gamma(H)$. The operation $\Omega_{2d+1}$ 
and $\Gamma_{2d+1}$ defined in Section 2 are graph functors.  

In the context of multiplicativity of graphs, the existence or non-existence of homomorphisms between two graphs (or digraphs) is more relevant than the structure of homomorphisms between them. 
Thus we consider the so called ``thin category" in which there is at most one morphism from one object to another. 

Given two digraph or graph functors $L$ and $R$, we say $L$ is a left adjoint of $R$ and $R$ is a right adjoint of $L$ if    $$L(G) \to H \Leftrightarrow G \to R(H).$$    

Assume   $L$ and $R$ form a pair of left and right adjoints. Then it is easy to verify that for any graphs (or digraphs) $G$ and $H$,  $R(G \times H) \leftrightarrow R(G) \times R(H)$. However, it is not necessarily true that $L(G \times H) \leftrightarrow L(G) \times L(H)$ for all graphs $G$ and $H$. If this is true, then  the right adjoint applied to a multiplicative graph results in a multiplicative graph.  

\begin{theorem}
    \label{thm-ft2018}
    If $L(G \times H) \leftrightarrow L(G) \times L(H)$ for all graphs (or digraphs)  $G$ and $H$, then $K$ being multiplicative implies that $R(K)$ is multiplicative. Moreover, if $L(R(G)) \leftrightarrow G$ for every graph (or digraph) $G$, then 
$R(K)$ being multiplicative implies that $K$ is multiplicative.  
\end{theorem}
  
For a positive integer $d$, let  $\Lambda_{2d+1}(G)$ be the graph obtained from $G$ by replacing each edge $e$ with a path of $2d+1$ edges. 
 It follows from results in \cite{Pultr1970,Tar2005,HT2010,FT2018} that
 \begin{enumerate}
     \item $\Gamma_{2d+1}(G \times H) \leftrightarrow \Gamma_{2d+1}(G) \times \Gamma_{2d+1}(H)$ and   $ \Omega_{2d+1}(G \times H) \leftrightarrow \Omega_{2d+1}(G) \times \Omega_{2d+1}(H)$,\\
     \item $\Gamma_{2d+1}(G) \to H \Leftrightarrow G \to \Omega_{2d+1}(H)$ and  $\Omega_{2d+1}(G) \to \Omega_{2d+1}(H) \Leftrightarrow G \to H$,\\
     \item $\Gamma_{2d+1}(\Omega_{2d+1}(G)) \leftrightarrow G \leftrightarrow \Gamma_{2d+1}(\Lambda_{2d+1}(G)).$
 \end{enumerate}
It follows from Theorem \ref{thm-ft2018} that  for any positive integer $d$, a graph  
 $K$ is multiplicative if and only if $\Omega_{2d+1}(K)$ is multiplicative.

 Note that $\Omega_{2d+1}(K_3) = \Omega_3(C_{2d+1}) = C_{3(2d+1)}$. Hence it follows from the multiplicativity of $K_3$  that for any positive integer $d$, $C_{3(2d+1)}$ is multiplicative, and from this it follows that $C_{2d+1}$ is multiplicative. 
 Hence all odd cycles are multiplicative. 
 This gives an alternate proof of Theorem \ref{thm-hhmn}.
Further applications of this graph functor lead to the proof of the multiplicativity of a family of {\em circular cliques},   which is related to the concept of circular chromatic number.
 
Recall that for positive integers $p \ge 2q$, a {  $(p,q)$-colouring} of a graph $G$ is a mapping $f: V(G) \to \{0,1,\ldots, p-1\}$ such that for every edge $xy$ of $G$, 
 $q \le |f(x)-f(y)| \le p-q.$ 
The {  circular chromatic number} of $G$ is  $\chi_c(G) = \inf\left\{ \frac pq: G \text{ admits a $(p,q)$-colouring} \right\}.$ 

Given positive integers $p \ge 2q$, the {\em circular  clique} $K_{p/q}$ is the graph with vertex  set $ \{0,1,\ldots, p-1\}$ and edge set $\{ij: q \le |i-j| \le p-q\}$. 

A $(p,q)$-colouring of a graph $G$ is equivalent to a homomorphism from $G$ to $K_{p/q}$. 
The circular clique $K_{2k+1/k}$ is the odd cycle $C_{2k+1}$. Thus the multiplicativity of $C_{2k+1}$ is equivalent to the statement that if $\chi_c(G), \chi_c(H) > 2+ 1/k$, then $\chi_c(G \times H) > 2+1/k$. 

As a strengthening of Hedetniemi's conjecture, this author conjectured in \cite{Zhu1992} that for any graphs $G$ and $H$, $\chi_c(G \times H) = \min\{\chi_c(G), \chi_c(H)\}$. This conjecture is equivalent to saying that the circular cliques $K_{p/q}$ are multiplicative.
This conjecture is stronger than 
  Hedetniemi's conjecture, and is hence false.  Nevertheless, Tardif \cite{Tar2005} proved  the following result.

  \begin{enumerate}
      \item If $p/q < 3$, then $\Gamma_3(K_{p/q}) \leftrightarrow K_{p/(3q-p)}$.
      \item If $p/q < 12/5$, then $\Omega_3(\Gamma_3(K_{p/q}) )\leftrightarrow K_{p/q}$.
  \end{enumerate}

It follows that for $p/q < 12/5$, $K_{p/q}$ is multiplicative if and only if $K_{p/(3q-p)}$ is multiplicative. Starting from the multiplicativity of $C_{2d+1}=K_{(2d+1)/d}$, one can infer that for a dense subset $S$ of rationals in $[2,4)$, $K_{p/q}$ is multiplicative for any $p/q \in S$. 
This implies that for any $p/q \in [2,4)$, 
if $\chi_c(G), \chi_c(H) > p/q$, then $\chi_c(G \times H) > p/q$. Thus Tardif \cite{Tar2005} proved the multiplicativity of all $K_{p/q}$ with $p/q < 4$.

\begin{theorem}\label{thm-kpq}
    If $p/q < 4$, then $K_{p/q}$ is multiplicative.
\end{theorem}

Theorem 
\ref{thm-kpq} is equivalent to the statement that if 
$\min\{\chi_c(G), \chi_c(H)\} < 4$, then
$\chi_c(G \times H) = \min\{\chi_c(G), \chi_c(H)\}$.

All the multiplicative graphs listed above have chromatic number at most $4$.
A natural question is whether there are multiplicative graphs of large chromatic number. 
This question was answered by Wrochna \cite{Wro2017}, who proved the following result that greatly enlarged the family of multiplicative graphs.

\begin{theorem}
    \label{thm-c4free}
    Any graph $G$ containing no 4-cycle is multiplicative.
\end{theorem}

In particular graphs of girth at least 5 are multiplicative. Hence  there are multiplicative graphs of arbitrarily large chromatic number. 
Theorem \ref{thm-c4free}   was further improved by Tardif and Wrochna in \cite{TW2019}. 

\begin{theorem}
    \label{thm-c4free+}
   If every edge of $G$ is contained in at most one copy of $C_4$, then $G$ is multiplicative. If $G$ has girth at least $13$, then $\Gamma_3(G)$ is multiplicative.  
\end{theorem}

The proofs of the multiplicativity of all the known multiplicative graphs $K$ have a feature in common and a seemingly much stronger result is proved: Assume $G$ and $H$ are connected graphs. Let $C$ and $D$ be odd cycles in $G$ and $H$ respectively. Then  
	$$(G \times D) \cup (H \times C) \to K \Rightarrow G \to K \text{ \ or \ } H \to K.$$

We call $K$ {\em strongly multiplicative} if the above statement is true.
It was conjectured in \cite{ES1985} that if $C$ and $D$ are $n$-chromatic subgraphs in $G$ and $H$ respectively, then
$(G \times D) \cup (C \times H) \to K_n$ implies that $G \to K_n$ or $H \to K_n$. This was proved for $n=3$ in \cite{ES1985}, and disproved for $n > 3$ in \cite{TZ2002}.

Multiplicativity of digraphs was first studied by Ne\v{s}et\v{r}il and Pultr in \cite{NP1978}, where the term ``productive" was used.   A family $\mathcal{C}$ of digraphs is defined to be  {\em productive} if it is closed under categorical product. 
For a given digraph $K$, denote by $\not\to K$ the family of digraphs $G$ such that $G \not\to K$. The productivity of classes of the form $\not\to K$, which is the same as the multiplicativity of $K$,  was studied in \cite{NP1978}.  

It  was proved in \cite{NP1978} that
    transitive tournaments, directed paths and directed cycles of prime length are multiplicative, and an oriented path is multiplicative if and only if it is homomorphically equivalent to a directed path.

 It was also observed in \cite{NP1978}  that a directed cycle whose length is not a prime power is not multiplicative, and  conjectured  that directed cycles of prime power length are multiplicative. This conjecture was proved by H\"{a}ggkvist,  Hell,  Miller, and   Neumann Lara  in \cite{HHMN1988} by using tools from homotopy theory. 
 
 \begin{theorem}
     \label{thm-dicycle}
     A directed cycle is multiplicative if and only if its length is a   prime power.
 \end{theorem}
 
 A combinatorial proof of Theorem \ref{thm-dicycle} was given by Zhou in   \cite{ZhouThesis}, and a simple combinatorial proof was given by this author in \cite{Zhu1992b}. 
 
 Multiplicativityy of oriented cycles is more complicated and studied in a few papers \cite{ZhouThesis, ZhuThesis, HZZ1994}. A complete characterization of multiplicative oriented cycles was given in \cite{HZZ1994}.  

Another type of digraph whose multiplicativity was studied is the class of acyclic local tournaments, which are acyclic digraphs $K$ such that for each vertex $v$, the in-neighbours of $v$ induce a tournament, and the out-neighbours of $v$ also induces a tournament. Transitive tournaments and directed paths are examples of acyclic local tournaments. A complete characterization of multiplicative acyclic local tournaments was given in \cite{ZZ1997}. 

To prove that a digraph $K$ is multiplicative, one usually needs a characterization of   digraphs that admit a homomorphism to $K$, in the form of homomorphism duality: $$ \to K = \mathcal{O} \not\to.$$
This means that a digraph $G$ admits a homomorphism to $K$ if and only if   there is no
homomorphism of any digraph in $\mathcal{O}$ to $G$.  For this duality to be useful, digraphs in $\mathcal{O}$ need to have a ``simple structure", and membership can be easily checked. 
For example, it was proved in \cite{NP1978} that if $K$ is a directed cycle of length $n$, then $\mathcal{O}$ consists of oriented cycles $C$ of length $\ell(C) \not\equiv 0 \pmod{n}$. Here, $\ell(C)$ is the number of forward arcs minus the number of backward arcs along a given transversal. This leads to a polynomial time algorithm for the $K$-colouring problem, and is   helpful in the proof of multiplicativity of $K$ when $K$ is   a directed cycle of prime power length $p^k$.  For that purpose, it suffices to show that if $C_1, C_2$ are oriented cycles of lengths $n_1, n_2 \not\equiv 0 \pmod{p^k}$, then there is an oriented cycle $C$ of length $n = {\rm lcm}(n_1,n_2)$ such that $C \to C_1 \times C_2$. 

Homomorphism duality is closely related to the complexity of the $H$-colouring problem, i.e., for a fixed graph or digraph $H$, decide whether an instance graph or digraph $G$ admits a homomorphism to $H$.
The complexity of the $H$-colouring problem for undirected graphs $H$ was solved by Hell and Ne\v{s}et\v{r}il \cite{HN1990}:

\begin{theorem}
    \label{thm-HN}
    The $H$-colouring problem is polynomial time solvable if $H$ is bipartite and NP-complete otherwise.
\end{theorem}

For digraphs $H$, the $H$-colouring problem is more complicated. It was proved by Feder and Vardi \cite{FV1999} that the complexity problem for a much more general class of problems, the {\em constraint satisfaction problems}  CSP, can be reduced to the complexity of $H$-colouring problems for digraphs $H$. They made the famous dichotomy conjecture:   every CSP is either polynomial time solvable or NP-complete. This is equivalent to  the CSP dichotomy conjecture restricted to digraphs, i.e.,  for any digraph $H$, the $H$-colouring problem is either polynomial time solvable or NP-complete.  The CSP has a wide range   of applications. The dichotomy conjecture had been a central problem in many different fields of mathematics and theoretical computer science for many years.  A wide variety of tools ranging from statistical physics to universal algebra  has been employed. The dichotomy conjecture was finally proved independently by Zhuk \cite{Zhuk2020} and by Bulatov \cite{Bul2017}  using algebraic methods. 

Multiplicativity of graphs and digraphs are also studied under the framework of lattice theory.   For graphs $G$ and $H$,  denote by $G+H$ the disjoint union of $G$ and $H$. Then we have $G, H \to G +H$.  
Let $\tilde{1}$ be the  graph with a single vertex  and a loop, and $\tilde{0}$ be the  graph with a single vertex and no edge.   The homomorphism relation $\to$ defines a Boolean lattice on $\mathcal{G}/ \leftrightarrow$, with $[G \times H]$ being the {\em meet} and $[G+H]$ being the {\em join} of $[G]$ and $[H]$, respectively  \cite{DS1996}, and with $[\tilde{0}]$ and $[\tilde{1}]$ as the $0$ and $1$ of the lattice, respectively.   

Note that $ G \to K$ if and only if $(G + K) \leftrightarrow K$. If $K$ is multiplicative, then $G \times H \to K$  implies that $G \to K$ or $H \to K$. Equivalently $(G+K) \leftrightarrow K$ or $(H+K) \leftrightarrow K$. Thus $K$ is multiplicative if and only if $[K]$ is {\em meet irreducible} in the lattice $\mathcal{G}/\leftrightarrow$, i.e., if $[K]$ is the meet of $[A]$ and $[B]$, then $[K] = [A]  $ or $[K] = [B]$.

Given a  graph $K$, $K^{\mathcal G}= \{[K^G]: G \in \mathcal{G}\}$ is also a Boolean lattice,  called an {\em exponential lattice} \cite{Tar2011}.  For a digraph $K$, $K^{\mathcal D}$ is defined similarly. The exponential lattice $K^{\mathcal G}$ is not a sub-lattice of $\mathcal {G}$. The lattice  $K_n^{\mathcal G}$ is closed under categorical product: $K^G \times K^H \leftrightarrow K^{G+H}$. However, it is not closed under addition: for $[G], [H] \in K^{\mathcal G}$, we usually do not have $[G + H] \in K^{\mathcal G}$. Nevertheless, 
$K^{\mathcal G}$ is a 
  Boolean lattice, where for $[G], [H] \in K^{\mathcal G}$, the join of $[G]$ and $[H]$ in $K^{\mathcal G}$ is $[K^{K^{G+H}}]$, which is generally larger than $G+H$ in the homomorphism order. 

  For any graph $K$, 
  $K^{\mathcal G}$ contains at least two elements: If $G$ is obtained from $K$ by adding a universal vertex, then $K^G \leftrightarrow K$. If $G \to K$, then 
$K^G \leftrightarrow \tilde{1}$. If $K$ is multiplicative, then $K^{\mathcal G}$ contains only these two elements. The converse is also true. Thus a graph $K$ is multiplicative if and only if $K^{\mathcal G}$ is a two-element lattice. Similarly,  a digraph $K$ is multiplicative if and only if $K^{\mathcal D}$ is a two-element lattice. 

As $K_n$ is not multiplicative for $n \ge 4$, 
$K_n^{\mathcal G}$ has more than two elements. Indeed, Theorem \ref{thm-Tardif2} implies that $K_n^{\mathcal{G}}$ is infinite. This is stronger than just saying that $K_n$ is non-multiplicative.
The following definition was given in \cite{TZ2002b}.

\begin{definition}
    \label{def-level}
    A graph $K$ is called $n$-composite if there exist graphs $H_1,H_2, \ldots, H_n$ such that $\prod_{i=1}^nH_i \to K$ and $\prod_{i \in I} H_i \not\to K$ for every proper subset $I$ of $\{1,2,\ldots, n\}$. The {\em level of non-multiplicativity } of $K$ is the largest integer $n$ such that $K$ is $n$-composite, if such an integer $n$ exists. Otherwise, $K$ is said to have infinite level of non-multiplicativity. 
\end{definition}
 
It was proved in \cite{TZ2002b} that   the countably infinite complete graph $K_{\aleph_0}$ has an infinite level of non-multiplicativity and that there exist Kneser graphs with arbitrarily high levels of non-multiplicativity. Theorem \ref{thm-Tardif2} implies that for $n \ge 4$, $K_n$ has an infinite level of  non-multiplicativity.

\section{The Poljak-R\"{o}dl function}

The Poljak-R\"{o}dl function is defined in \cite{PR1981} as follows:
$$f(n) = \min\{\chi(G \times H): \chi(G), \chi(H) \ge n\}.$$
The statement $H(n)$ is equivalent to  $f(n+1) =n+1$. As $H(n)$ fails for $n \ge 4$, it follows that $f(n) \le n-1$ for $n \ge 5$. A natural question is how large can be the difference $n-f(n)$. 
Using Shitov's result,
Tardif and Zhu \cite{TZ2019} proved that  $f(n) \le n - \left(\log n\right)^{1/4-o(1)}$. So the difference 
 $n-f(n)$ can be arbitrarily large.
Tardif and Zhu \cite{TZ2019} asked the  question   whether the ratio $f(n)/n$ can be bounded away from $1$.  He and Wigderson \cite{HW2021} answered this question   in affirmative by showing that $f(n) \le (1 - 10^{-9})n + o(n)$. Then it was proved in \cite{Zhu2020}  that  $f(n) \le  n/2 + o(n)$, and this result was further refined by Tardif in \cite{Tar2023}, where the following result was proved.

\begin{theorem} 
	\label{thm-poljakrodl} 
	  $f(m) \le \lceil m/2 \rceil +3$.
\end{theorem}

The proof of Theorem \ref{thm-poljakrodl} is quite simple.
Let $m=2n-6$,  $w$ be a sufficiently large  odd integer,  and let $G= \Omega_w(K_m)$. We shall show that $\chi(K_n^G) \ge m$. 
 
Assume $\chi(K_n^G) = k$ and $\Psi$ is a $k$-colouring of $K_n^G$. We may assume the constant mapping ${\rm const}_i$ is coloured by colour $i$. Thus, for any mapping $\phi$, $\Psi(\phi) \in Im(\phi) \cup \{n+1, n+2, \ldots, k\}$.

By Theorem \ref{thm-Tardif2}, $\chi(K_4^G) > 4$.
For a 4-subset $A$ of $[n]$, let  $H_A$ be the subgraph of $K_n^G$ induced by those mappings $\phi$ with   $Im(\phi) \subseteq A$. Then $H_A$ is isomorphic to $K_4^G$, and hence is not $4$-colourable. So $\Psi(H_A) \cap \{n+1,n+2, \ldots, k\} \ne \emptyset$.  Let $\tau(A) $ be an arbitrary colour in $\Psi(H_A) \cap \{n+1,n+2, \ldots, k\}$.  
	If $A, B$ are $4$-subsets of $[n]$ and $A \cap B = \emptyset$, then every vertex in $H_A$ is adjacent to every vertex in $H_B$. Hence $\Psi(H_A) \cap \Psi(H_B) = \emptyset$. In particular, $\tau(A) \ne \tau(B)$. Thus $\tau$ is a proper colouring of the Kneser graph $K(n, 4)$.
 By Lov\'{a}sz's Theorem, $\chi(K(n, 4)) = n-6$. Hence $k  \ge m= 2n-6$.  	

Now    $G= \Omega_w(K_m)$ and $K_n^G$ both have chromatic number at least $m$, and their product is $(\lceil m/2 \rceil +3)$-colourable. Hence 
$f(m) \le \lceil m/2 \rceil +3$.    
This completes the proof of Theorem \ref{thm-poljakrodl}.

On the other hand, the only known lower bound for $f(n)$ is $f(n) \ge 4$ for $n \ge 4$, a result that has not been improved in three decades. Notably, it remains unknown whether $f(n)$ is bounded by a constant. Nevertheless,  the following result was proved in \cite{Pol1991}:

\begin{theorem}
    \label{thm-poljak}
    $f(n)$ is either unbounded or bounded by $9$. 
\end{theorem}

Theorem \ref{thm-poljak}  was proved by considering   the product of digraphs.  
The chromatic number of a digraph $G$ is defined as $\chi(\underline{G})$, where $\underline{G}$ is the underlying  graph of $G$,  obtained by replacing each arc $(x,y)$ with an edge $xy$. We view the complete graph $K_n$ as a symmetric digraph, in which $(x,y)$ is an arc for every pair of distinct vertices $x$ and $y$. Then an $n$-colouring of a digraph $G$ is a homomorphism from $G$ to $K_n$.

Given a digraph $G$, let $G^{-1}$ be the digraph obtained from $G$ by reversing the direction of all its arcs. It is easy to see that 
for any digraphs $G$ and $H$, $$\underline{G} \times \underline{H} = (\underline{G \times H})
\cup (\underline{G \times H^{-1}}).$$	
Hence $$\chi (\underline{G} \times \underline{H}) \le \chi (G \times H) \times \chi (G \times H^{-1}).$$
Let 
\begin{eqnarray*}
	g(n) &=& \min\{\chi(G \times H): \text{ $G$ and $H$ are digraphs with } \chi(G),  \chi(H) \ge n\},\\
	h(n) &=& \min\{\max\{\chi(G \times H), \chi(G \times H^{-1})\}: \text{ $G$ and $H$ are digraphs with } \chi(G),  \chi(H) \ge n\}.
\end{eqnarray*}	
Since $E(\underline{G}\times \underline{H}) = E(\underline{G \times H}) \cup E(\underline{G \times H^{-1}})$, we have  
$$g(n) \le h(n) \le  f(n) \le h(n)^2.$$ 

Thus $f(n)$ is bounded by a constant if and only if $h(n)$ is bounded by a constant. 
Tardif and Wehlau \cite{TW2006} further proved that $f(n)$ is bounded by a constant if and only if $g(n)$ is bounded by a constant.

  Poljak and R\"{o}dl   \cite{PR1981} proved that if $g(n)$  (respectively $h(n)$) is bounded by a constant, then the smallest such constant is at most $4$.  
  This result was further improved  independently by Schmerl    (unpublished),  Poljak \cite{Pol1991}, and this author (unpublished), who proved that if $g(n)$  (respectively $h(n)$) is bounded by a constant, then the smallest such constant is at most $3$.
 Consequently, if $f(n)$ is bounded by a constant,  then the smallest  such constant is at most $9$.

If $K_n^{\mathcal G}$   contains an atom, then there is an integer $N$ such that $\chi(G), \chi(H) > N$ implies that $\chi(G \times H) > n$.  In particular, if $K_9^{\mathcal G}$ has an atom, then there is an integer $N$ such that 
$\chi(G), \chi(H) \ge N$ implies that $\chi(G \times H) \ge 10$, i.e., for the Poljak-R\"{o}dl defined $f(n)$ in Section 4, $f(N) \ge 10$. Hence by Theorem \ref{thm-poljak}, $f(n)$ is unbounded. 

Similarly, if $K_3^{\mathcal D}$ has an atom, then the function $g(n)$ defined in Section 4 is unbounded, and hence $f(n)$ is also unbounded. However, it was proved by Tardif \cite{Tar2011} that for $n \ge 3$, $K_3^{\mathcal D}$ does not have an atom, and Theorem \ref{thm-Tardif2} implies that for $n \ge 4$, $K_n^{\mathcal{G}}$ is infinite. On the other hand, this does not mean that $g(n)$ or $f(n)$ is bounded.



 \section{Open questions}
 
Although Hedetniemi's conjecture has been refuted,  the chromatic number of the product of two graphs   with given chromatic numbers remains  an interesting and challenging problem. One  natural question is whether the Poljak-R\"{o}dl function goes to infinity. The statement that $\lim_{n \to \infty}f(n) = \infty$ is referred as the weak Hedetniemi's conjecture. There is not much support for either the weak Hedetniemi's conjecture or its negation. 

If the weak Hedetniemi's conjecture is out of reach, then one naturally considers more restricted questions. Can the lower bound and upper bound for $f(n)$ be improved? 

\begin{question}
    \label{q1}
      Is 
    there an integer $n$ such that $f(n) \ge 5$? 
\end{question}

By Theorem \ref{thm-Tardif}, for $n \ge 5$, $f(n) \le n-1$. As $f(n) \ge 4$ for $n \ge 4$, we have  $f(5)=4$ and $f(6) =4$ or $5$. Can we determine $f(6)$?

\begin{question}
    \label{q2}
   Can the upper bound $f(m) \le \lceil m/2 \rceil +3$ be improved? 
   Is it true that $\lim \sup  f(n)/n = 0$?  
\end{question}

By Theorem \ref{thm-poljakrodl}, $f(m) \le \lceil m/2 \rceil +3$. For the proof of this result, let $m=2n-6$ and $w$ be a sufficiently large  odd integer and let $G= \Omega_w(K_m)$. It is proved that that $\chi(K_n^G) \ge m$. The proof actually shows that for any $m' \ge m$,
there is an odd integer $w$ such that for $G=\Omega_w(K_{m'})$,   $\chi(K_n^G) \ge m$. 
To improve the upper bound, it suffices  to prove that  for some $m' > m$, for $G = \Omega_w(K_{m'})$, $\chi(K_n^G) \ge m'$. In the proof of Theorem \ref{thm-poljakrodl}, a subgraph $H$ of $K_n^G$ consisting of  mappings $\phi \in K_n^G$ with $|Im(\phi)| \le 4$ are shown to have $\chi(H) \ge m$. 
It is likely that the whole graph  has   chromatic number much larger than $m$, or even the subgraph $H$ described above already has chromatic number much larger than $m$. By Theorem \ref{thm-Tardif2}, for any integer $m$, if $w$ is a large enough odd integer and $G=\Omega_w(K_m)$,
$\chi(K_4^G) \ge 5$.  It is unknown if $\chi(K_4^G)$ can be much larger than $5$ (if $w$ is large enough). 

Similar question for the function $$g(n) = \min\{\chi(G \times H): \text{ $G$ and $H$ are digraphs with $\chi(G) = \chi(H) = n$}\}$$   is equally interesting. However, this question is less studied.  It was already known   that $g(n) < n$ for $n \ge 3$ when the function $g(n)$ was defined in \cite{PR1981}. In 2004,   Tardif \cite{Tar2004} proved that there are $n$-tournaments $S$ and $T$ such that $\chi(S \times T) \le 2n/3$, which implies that $g(n) \le 2n/3$.  
As $g(n) \le f(n)$, the upper bound $f(n) \le \lceil n/2 \rceil +3$ is also an upper bound for $g(n)$.  It is natural that one expects a better upper bound for $g(n)$, however, currently $\lceil n/2 \rceil +3$ is the best known upper bound for $g(n)$. 
  
\begin{question}
    \label{q3} 
   Is it true that $\lim \sup  g(n)/n = 0$?  
\end{question}

Although Hedetniemi's conjecture is not true, the problem that for which graphs $G$ and $H$, $\chi(G \times H) = \min\{\chi(G), \chi(H)\}$ is still an interesting problem. We have some counterexamples to Hedetniemi's conjecture, however, it is not clear what makes a pair of graphs a counterexample to $H(n)$ or what makes a graph $G$ a counterexample to $H'(n)$? In all the known counterexamples $G$ to $H'(n)$, one finds a non $n$-colourable subgraph $H$ of $K_n^G$ that are induced by vertices in the neighbourhood of the $n$-clique induced by the constant mappings. Now $H$ is also a counterexample to $H'(n)$, and $G$ is a non-$n$-colourable subgraph of $K_n^H$. Then the connected component of $K_n^H$   containing  $G$ as a subgraph does not contain an $n$-clique. For otherwise, we would have two connected graphs $G_1,G_2$ with $\chi(G_i) > n$, each containing an $n$-clique and $\chi(G_1 \times G_2) =n$, which contradicts a result proved in \cite{DSW1985}.  Therefore $G$ and $H$ are both counterexamples to $H'(n)$, but seems to have quite different structures. Does this suggest that counterexamples to $H'(n)$ are abundant? 

\begin{question}
    \label{q4}
What is the probability that a random graph $G$ with $\chi(G) > n$ of large girth being a counterexample to $H'(n)$?
\end{question}

Tardif (personal communication) asked the following question.

\begin{question}
    \label{q5}
    For an odd cycle $C_{2d+1}$ for $d \ge 2$, are there  integers $m, n$ such that the lexicographic product $C[K_m]$ is a counterexample to $H'(n)$?  
\end{question}
 
 It was proved in \cite{DSW1985} that if $\chi(G) > n$ and each vertex of $G$ is contained in an $n$-clique, then $\chi(K_n^G) =n$. Each vertex of $C[K_m]$ is contained in a $2m$-clique. On the other hand, $\chi(C_{2d+1}[K_m]) = 2m+ \lceil m/d \rceil$. So if such integers $m,n$ exist, then $2m < n <   2m+ \lceil m/d \rceil$.

The chromatic number of the product of tournaments was studied by Tardif. It was proved in  \cite{Tar2004} that $ \min \{ \chi (S \times T): \text{ $S$ and $T$ are $n$-tournaments} \}$   is asymptotically equal to $\lambda n$, where $1/2 \le \lambda \le 2/3$. 

\begin{question}
    \label{q6} What is  $\lim_{n \to \infty} \min \{ \chi (S \times T): \text{ $S$ and $T$ are $n$-tournaments} \}$ ? 
\end{question}

\begin{question}
    \label{q7} 
    For which homomorphism-monotone graph invariants $\rho$, the Hedetniemi-type equality 
    $\rho(G \times H) = \min \{\rho(G), \rho(H)\}$ holds? 
\end{question}

  Elements in the asymptotic spectrum  $Y(\mathcal{G})$ of the family $\mathcal{G}$ of all graphs have nice properties with respect to the OR-product and the join operations.   It was remarked in \cite{Zui2019} that this family of graph invariants has infinitely many different
elements. The Hedetniemi type equality are studied for a few elements in this family, and the
corresponding problems are challenging and interesting. This could be the case for some other graph invariants in this family.

 \section*{Acknowledgement}

I thank Claude Tardif  for many valuable comments. I thank  Tao Wang  and anonymous referees  for  careful readings of the manuscript and many helpful comments that improved the presentation of the paper.
 
 \bibliographystyle{abbrv}
\bibliography{ref}

\begin{thebibliography}{10}

\bibitem{Alon1998}
N.~Alon.
\newblock The {S}hannon capacity of a union.
\newblock {\em Combinatorica}, 18(3):301--310, 1998.

\bibitem{BS2005}
S.~Baum and M.~Stiebitz.
\newblock Coloring of graphs without short odd paths between vertices of the
  same color class.
\newblock {\em Unpublished manuscript}, 2005.

\bibitem{Bul2017}
A.~A. Bulatov.
\newblock A dichotomy theorem for nonuniform {CSP}s.
\newblock In {\em 58th {A}nnual {IEEE} {S}ymposium on {F}oundations of
  {C}omputer {S}cience---{FOCS} 2017}, pages 319--330. IEEE Computer Soc., Los
  Alamitos, CA, 2017.

\bibitem{BEL1976}
S.~A. Burr, P.~Erd\H{o}s, and L.~Lov\'{a}sz.
\newblock On graphs of {R}amsey type.
\newblock {\em Ars Combin.}, 1(1):167--190, 1976.

\bibitem{Cso2007}
P.~Csorba.
\newblock Homotopy types of box complexes.
\newblock {\em Combinatorica}, 27(6):669--682, 2007.

\bibitem{Cso2008}
P.~Csorba.
\newblock On the simple {$\Bbb Z_2$}-homotopy types of graph complexes and
  their simple {$\Bbb Z_2$}-universality.
\newblock {\em Canad. Math. Bull.}, 51(4):535--544, 2008.

\bibitem{DRV2023}
H.~R. Daneshpajouh, R.~Karasev, and A.~Volovikov.
\newblock Hedetniemi's conjecture from the topological viewpoint.
\newblock {\em J. Combin. Theory Ser. A}, 195:Paper No. 105721, 24, 2023.

\bibitem{DSW1985}
D.~Duffus, B.~Sands, and R.~E. Woodrow.
\newblock On the chromatic number of the product of graphs.
\newblock {\em J. Graph Theory}, 9(4):487--495, 1985.

\bibitem{DS1996}
D.~Duffus and N.~Sauer.
\newblock Lattices arising in categorial investigations of {H}edetniemi's
  conjecture.
\newblock {\em Discrete Math.}, 152(1-3):125--139, 1996.

\bibitem{ES1985}
M.~El-Zahar and N.~W. Sauer.
\newblock The chromatic number of the product of two {$4$}-chromatic graphs is
  {$4$}.
\newblock {\em Combinatorica}, 5(2):121--126, 1985.

\bibitem{Erd1959}
P.~Erd\H{o}s.
\newblock Graph theory and probability.
\newblock {\em Canadian J. Math.}, 11:34--38, 1959.

\bibitem{FV1999}
T.~Feder and M.~Y. Vardi.
\newblock The computational structure of monotone monadic {SNP} and constraint
  satisfaction: a study through {D}atalog and group theory.
\newblock {\em SIAM J. Comput.}, 28(1):57--104, 1999.

\bibitem{FT2018}
J.~Foniok and C.~Tardif.
\newblock Hedetniemi's conjecture and adjoint functors in thin categories.
\newblock {\em Appl. Categ. Structures}, 26(1):113--128, 2018.

\bibitem{GRRSV2020}
C.~Godsil, D.~E. Roberson, B.~Rooney, R.~\v{S}\'{a}mal, and A.~Varvitsiotis.
\newblock Vector coloring the categorical product of graphs.
\newblock {\em Math. Program.}, 182(1-2):275--314, 2020.

\bibitem{GRSS2016}
C.~Godsil, D.~E. Roberson, R.~\v{S}\'{a}mal, and S.~Severini.
\newblock Sabidussi versus {H}edetniemi for three variations of the chromatic
  number.
\newblock {\em Combinatorica}, 36(4):395--415, 2016.

\bibitem{GS2022}
A.~Gujgiczer and G.~Simonyi.
\newblock On multichromatic numbers of widely colorable graphs.
\newblock {\em J. Graph Theory}, 100(2):346--361, 2022.

\bibitem{GJS2004}
A.~Gy\'{a}rf\'{a}s, T.~Jensen, and M.~Stiebitz.
\newblock On graphs with strongly independent color-classes.
\newblock {\em J. Graph Theory}, 46(1):1--14, 2004.

\bibitem{Hae1979}
W.~Haemers.
\newblock On some problems of {L}ov\'asz concerning the {S}hannon capacity of a
  graph.
\newblock {\em IEEE Trans. Inform. Theory}, 25(2):231--232, 1979.

\bibitem{HHMN1988}
R.~H\"{a}ggkvist, P.~Hell, D.~J. Miller, and V.~Neumann~Lara.
\newblock On multiplicative graphs and the product conjecture.
\newblock {\em Combinatorica}, 8(1):63--74, 1988.

\bibitem{Haj2009}
H.~Hajiabolhassan.
\newblock On colorings of graph powers.
\newblock {\em Discrete Math.}, 309(13):4299--4305, 2009.

\bibitem{HT2010}
H.~Hajiabolhassan and A.~Taherkhani.
\newblock Graph powers and graph homomorphisms.
\newblock {\em Electron. J. Combin.}, 17(1):Research Paper 17, 16, 2010.

\bibitem{HW2021}
X.~He and Y.~Wigderson.
\newblock Hedetniemi's conjecture is asymptotically false.
\newblock {\em J. Combin. Theory Ser. B}, 146:485--494, 2021.

\bibitem{Hed1966}
S.~T. Hedetniemi.
\newblock {\em Homomoprhisms of graphs and automata}.
\newblock ProQuest LLC, Ann Arbor, MI, 1966.
\newblock Thesis (Ph.D.)--University of Michigan.

\bibitem{Hell1979}
P.~Hell.
\newblock An introduction to the category of graphs.
\newblock In {\em Topics in graph theory ({N}ew {Y}ork, 1977)}, volume 328 of
  {\em Ann. New York Acad. Sci.}, pages 120--136. New York Acad. Sci., New
  York, 1979.

\bibitem{HN1990}
P.~Hell and J.~Ne\v{s}et\v{r}il.
\newblock On the complexity of {$H$}-coloring.
\newblock {\em J. Combin. Theory Ser. B}, 48(1):92--110, 1990.

\bibitem{HNbook}
P.~Hell and J.~Ne\v{s}et\v{r}il.
\newblock {\em Graphs and homomorphisms}.
\newblock Oxford University Press, 2004.

\bibitem{HZZ1994}
P.~Hell, H.~Zhou, and X.~Zhu.
\newblock Multiplicativity of oriented cycles.
\newblock {\em J. Combin. Theory Ser. B}, 60(2):239--253, 1994.

\bibitem{KMS1998}
D.~Karger, R.~Motwani, and M.~Sudan.
\newblock Approximate graph coloring by semidefinite programming.
\newblock {\em J. ACM}, 45(2):246--265, 1998.

\bibitem{Kla1996}
S.~Klav\v{z}ar.
\newblock Coloring graph products---a survey.
\newblock {\em Discrete Math.}, 155(1-3):135--145, 1996.
\newblock Combinatorics (Acireale, 1992).

\bibitem{Lov1978}
L.~Lov\'{a}sz.
\newblock Kneser's conjecture, chromatic number, and homotopy.
\newblock {\em J. Combin. Theory Ser. A}, 25(3):319--324, 1978.

\bibitem{Lov1979}
L.~Lov\'asz.
\newblock On the {S}hannon capacity of a graph.
\newblock {\em IEEE Trans. Inform. Theory}, 25(1):1--7, 1979.

\bibitem{MR708798}
L.~Lov\'asz.
\newblock Self-dual polytopes and the chromatic number of distance graphs on
  the sphere.
\newblock {\em Acta Sci. Math. (Szeged)}, 45(1-4):317--323, 1983.

\bibitem{MR1988723}
J.~Matou\v{s}ek.
\newblock {\em Using the {B}orsuk-{U}lam theorem}.
\newblock Universitext. Springer-Verlag, Berlin, 2003.
\newblock Lectures on topological methods in combinatorics and geometry,
  Written in cooperation with Anders Bj\"orner and G\"unter M. Ziegler.

\bibitem{Matsushita}
T.~Matsushita.
\newblock {$\Bbb Z_2$}-indices and {H}edetniemi's conjecture.
\newblock {\em Discrete Comput. Geom.}, 62(3):662--673, 2019.

\bibitem{NP1978}
J.~Ne\v{s}et\v{r}il and A.~Pultr.
\newblock On classes of relations and graphs determined by subobjects and
  factorobjects.
\newblock {\em Discrete Math.}, 22(3):287--300, 1978.

\bibitem{NZ2004}
J.~Ne\v{s}et\v{r}il and X.~Zhu.
\newblock On sparse graphs with given colorings and homomorphisms.
\newblock {\em J. Combin. Theory Ser. B}, 90(1):161--172, 2004.
\newblock Dedicated to Adrian Bondy and U. S. R. Murty.

\bibitem{Pol1991}
S.~Poljak.
\newblock Coloring digraphs by iterated antichains.
\newblock {\em Comment. Math. Univ. Carolin.}, 32(2):209--212, 1991.

\bibitem{PR1981}
S.~Poljak and V.~R\"{o}dl.
\newblock On the arc-chromatic number of a digraph.
\newblock {\em J. Combin. Theory Ser. B}, 31(2):190--198, 1981.

\bibitem{Pultr1970}
A.~Pultr.
\newblock The right adjoints into the categories of relational systems.
\newblock In {\em Reports of the {M}idwest {C}ategory {S}eminar, {IV}}, volume
  Vol. 137 of {\em Lecture Notes in Math.}, pages 100--113. Springer,
  Berlin-New York, 1970.

\bibitem{Sau2001}
N.~Sauer.
\newblock Hedetniemi's conjecture---a survey.
\newblock {\em Discrete Math.}, 229(1-3):261--292, 2001.

\bibitem{Schrijver1979}
A.~Schrijver.
\newblock A comparison of the {D}elsarte and {L}ov\'asz bounds.
\newblock {\em IEEE Trans. Inform. Theory}, 25(4):425--429, 1979.

\bibitem{Shi2019}
Y.~Shitov.
\newblock Counterexamples to {H}edetniemi's conjecture.
\newblock {\em Ann. of Math. (2)}, 190(2):663--667, 2019.

\bibitem{STW2017}
G.~Simons, C.~Tardif, and D.~Wehlau.
\newblock Generalised mycielski graphs, signature systems, and bounds on
  chromatic numbers.
\newblock {\em J. Combinatorial Theory Ser. B}, 122(2):776--793, 2017.

\bibitem{Sim2021}
G.~Simonyi.
\newblock Shannon capacity and the categorical product.
\newblock {\em Electron. J. Combin.}, 28(1):Paper No. 1.51, 23, 2021.

\bibitem{ST2006}
G.~Simonyi and G.~Tardos.
\newblock Local chromatic number, {K}y {F}an's theorem and circular colorings.
\newblock {\em Combinatorica}, 26(5):587--626, 2006.

\bibitem{SZ2010}
G.~Simonyi and A.~Zsb\'{a}n.
\newblock On topological relaxations of chromatic conjectures.
\newblock {\em European J. Combin.}, 31(8):2110--2119, 2010.

\bibitem{Tar2004}
C.~Tardif.
\newblock Chromatic numbers of products of tournaments: fractional aspects of
  {H}edetniemi's conjecture.
\newblock In {\em Graphs, morphisms and statistical physics}, volume~63 of {\em
  DIMACS Ser. Discrete Math. Theoret. Comput. Sci.}, pages 171--175. Amer.
  Math. Soc., Providence, RI, 2004.

\bibitem{Tar2005a}
C.~Tardif.
\newblock The fractional chromatic number of the categorical product of graphs.
\newblock {\em Combinatorica}, 25(5):625--632, 2005.

\bibitem{Tar2005}
C.~Tardif.
\newblock Multiplicative graphs and semi-lattice endomorphisms in the category
  of graphs.
\newblock {\em J. Combin. Theory Ser. B}, 95(2):338--345, 2005.

\bibitem{Tar2008}
C.~Tardif.
\newblock Hedetniemi's conjecture, 40 years later.
\newblock {\em Graph Theory Notes N. Y.}, 54:46--57, 2008.

\bibitem{Tar2011}
C.~Tardif.
\newblock Hedetniemi's conjecture and dense {B}oolean lattices.
\newblock {\em Order}, 28(2):181--191, 2011.

\bibitem{Tar2022c}
C.~Tardif.
\newblock The chromatic number of the product of 14-chromatic graphs can be 13.
\newblock {\em Combinatorica}, 42(2):301--308, 2022.

\bibitem{Tar2022}
C.~Tardif.
\newblock Counterexamples to {H}edetniemi's conjecture and infinite {B}oolean
  lattices.
\newblock {\em Comment. Math. Univ. Carolin.}, 63(3):315--327, 2022.

\bibitem{Tar2023}
C.~Tardif.
\newblock The chromatic number of the product of 5-chromatic graphs can be 4.
\newblock {\em Combinatorica}, 43(6):1067--1073, 2023.

\bibitem{Tar2024}
C.~Tardif.
\newblock Chromatic ramsey numbers of generalized mycielski graphs.
\newblock {\em Discussiones Mathematicae Graph Theory}, 44, 2024.

\bibitem{TW2006}
C.~Tardif and D.~Wehlau.
\newblock Chromatic numbers of products of graphs: the directed and undirected
  versions of the {P}oljak-{R}\"{o}dl function.
\newblock {\em J. Graph Theory}, 51(1):33--36, 2006.

\bibitem{TW2019}
C.~Tardif and M.~Wrochna.
\newblock Hedetniemi's conjecture and strongly multiplicative graphs.
\newblock {\em SIAM J. Discrete Math.}, 33(4):2218--2250, 2019.

\bibitem{TZ2002b}
C.~Tardif and X.~Zhu.
\newblock The level of nonmultiplicativity of graphs.
\newblock {\em Discrete Math.}, 244(1-3):461--471, 2002.

\bibitem{TZ2002}
C.~Tardif and X.~Zhu.
\newblock On {H}edetniemi's conjecture and the colour template scheme.
\newblock {\em Discrete Math.}, 253(1-3):77--85, 2002.

\bibitem{TZ2019}
C.~Tardif and X.~Zhu.
\newblock A note on {H}edetniemi's conjecture, {S}tahl's conjecture and the
  {P}oljak-{R}\"{o}dl function.
\newblock {\em Electron. J. Combin.}, 26(4):Paper No. 4.32, 5, 2019.

\bibitem{Vin1988}
A.~Vince.
\newblock Star chromatic number.
\newblock {\em J. Graph Theory}, 12(4):551--559, 1988.

\bibitem{Ziv2005}
R.~T. \v{Z}ivaljevi\'{c}.
\newblock W{I}-posets, graph complexes and {${\Bbb Z}_2$}-equivalences.
\newblock {\em J. Combin. Theory Ser. A}, 111(2):204--223, 2005.

\bibitem{Wro2017}
M.~Wrochna.
\newblock Square-free graphs are multiplicative.
\newblock {\em J. Combin. Theory Ser. B}, 122:479--507, 2017.

\bibitem{Wro2019}
M.~Wrochna.
\newblock On inverse powers of graphs and topological implications of
  {H}edetniemi's conjecture.
\newblock {\em J. Combin. Theory Ser. B}, 139:267--295, 2019.

\bibitem{Wro2020}
M.~Wrochna.
\newblock Smaller counterexamples to {H}edetniemi's conjecture.
\newblock {\em arXiv:2012.13558}, 2020.

\bibitem{Zhang2012}
H.~Zhang.
\newblock Independent sets in direct products of vertex-transitive graphs.
\newblock {\em J. Combin. Theory Ser. B}, 102(3):832--838, 2012.

\bibitem{ZhouThesis}
H.~Zhou.
\newblock Homomorphism properties of graph products.
\newblock {\em Ph.D. thesis, Simon Fraser University, Burnaby, B.C}, 1988.

\bibitem{ZZ1997}
H.~Zhou and X.~Zhu.
\newblock Multiplicativity of acyclic local tournaments.
\newblock {\em Combinatorica}, 17(1):135--145, 1997.

\bibitem{ZhuThesis}
X.~Zhu.
\newblock Multiplicative structures.
\newblock {\em Ph.D. thesis, The University of Calgary}, 1990.

\bibitem{Zhu1992b}
X.~Zhu.
\newblock A simple proof of the multiplicativity of directed cycles of prime
  power length.
\newblock {\em Discrete Appl. Math.}, 36(3):313--316, 1992.

\bibitem{Zhu1992}
X.~Zhu.
\newblock Star chromatic numbers and products of graphs.
\newblock {\em J. Graph Theory}, 16(6):557--569, 1992.

\bibitem{Zhu1998}
X.~Zhu.
\newblock A survey on {H}edetniemi's conjecture.
\newblock {\em Taiwanese J. Math.}, 2(1):1--24, 1998.

\bibitem{ZhuSurvey2001}
X.~Zhu.
\newblock Circular chromatic number: a survey.
\newblock {\em Discrete Math.}, 229(1-3):371--410, 2001.

\bibitem{Zhu2002}
X.~Zhu.
\newblock The fractional chromatic number of the direct product of graphs.
\newblock {\em Glasg. Math. J.}, 44(1):103--115, 2002.

\bibitem{Zhu2011}
X.~Zhu.
\newblock The fractional version of {H}edetniemi's conjecture is true.
\newblock {\em European J. Combin.}, 32(7):1168--1175, 2011.

\bibitem{Zhu2020}
X.~Zhu.
\newblock A note on the {P}oljak-{R}\"{o}dl function.
\newblock {\em Electron. J. Combin.}, 27(3):Paper No. 3.2, 4, 2020.

\bibitem{Zhu2021}
X.~Zhu.
\newblock Relatively small counterexamples to {H}edetniemi's conjecture.
\newblock {\em J. Combin. Theory Ser. B}, 146:141--150, 2021.

\bibitem{Zhuk2020}
D.~Zhuk.
\newblock A proof of the {CSP} dichotomy conjecture.
\newblock {\em J. ACM}, 67(5):Art. 30, 78, 2020.

\bibitem{Zui2019}
J.~Zuiddam.
\newblock The asymptotic spectrum of graphs and the {S}hannon capacity.
\newblock {\em Combinatorica}, 39(5):1173--1184, 2019.

\end{thebibliography}

\end{document}